\documentclass[12pt]{amsart}

\def\Xint#1{\mathchoice
{\XXint\displaystyle\textstyle{#1}}%
{\XXint\textstyle\scriptstyle{#1}}%
{\XXint\scriptstyle\scriptscriptstyle{#1}}%
{\XXint\scriptscriptstyle\scriptscriptstyle{#1}}%
\!\int}
\def\XXint#1#2#3{{\setbox0=\hbox{$#1{#2#3}{\int}$ }
\vcenter{\hbox{$#2#3$ }}\kern-.6\wd0}}

\def\dashint{\Xint-}

\usepackage{pgf,tikz}
\usepackage{mathrsfs}
\usepackage[font={footnotesize,it}]{caption}
\usepackage{float}

\usepackage[hypertexnames=false, colorlinks, citecolor=red, linkcolor=red, pagebackref]{hyperref}
\usepackage{geometry,amsmath,amssymb,enumerate,bbm}
\usepackage[american]{babel}
\usepackage{comment}
\usepackage{centernot}
\usepackage{graphicx}
\usepackage{psfrag}
\usepackage{amsmath} 
\usepackage{amsfonts}
\usepackage{amssymb}
\usepackage{amsthm}
\usepackage[normalem]{ulem}

\allowdisplaybreaks

\newcommand{\norm}[1]{\Vert #1 \Vert}
\newcommand{\assign}{:=}
\newcommand{\backassign}{=:}

\newcommand{\tmem}[1]{{\em #1\/}}
\newcommand{\tmop}[1]{\ensuremath{\operatorname{#1}}}

\newenvironment{enumerateroman}{\begin{enumerate}[i.] }{\end{enumerate}}
\newtheorem{theorem}{Theorem}
\newtheorem{lemma}[theorem]{Lemma}

\newtheorem{problem}{Problem}
\newtheorem{proposition}[theorem]{Proposition}

\newtheorem{example}[theorem]{Example}

\theoremstyle{definition}

\renewenvironment{proof}[1][Proof]{\textbf{#1.} }{\ \rule{0.5em}{0.5em}}

\newcommand{\III}{{\mathcal{I}}}

\newcommand{\RR}{{\mathbb{R}}}

\newcommand{\CCap}{\text{Cap}}

\newcommand{\pp}{{p^*}}

\DeclareMathOperator{\Capac}{Cap}

\begin{document}

\title{Two-weight dyadic Hardy's inequalities}

\author[]{N. Arcozzi, N. Chalmoukis, M. Levi, P. Mozolyako}
 \address{Dipartimento di Matematica, Alma Mater Studorium Universit\`a di
  Bologna, Piazza di Porta San Donato 5, 40126 Bologna, Italy}
  \email{{\tt nicola.arcozzi@unibo.it}}
  \thanks{The first three authors are members of INdAM}
 \address{Faculty of Mathematics and Computer Science, Saarland University, Saarbrücken Campus
66123 Saarbrücken, Germany}
  \email{{\tt chalmoukis@math.uni-sb.de}}
  \thanks{N. Chalmoukis was supported by the Alexander von Humboldt foundation, Humboldt Research Fellowship for postdoctoral and experienced researchers}
   \address{MaLGa Center - DIBRIS  - Università di Genova, Via Dodecaneso, 35, 16146 Genova, Italy.}
  \email{{\tt m.l.matteolevi@gmail.com}}
  \thanks{ M. Levi was supported by the Project “Harmonic Analysis on Continuous and Discrete Structures” (bando Trapezio Compagnia di San Paolo CUP E13C21000270007) and by the Project GNAMPA 2022 ``Generalized Laplacians on Continuous and Discrete Structures" (CUP E55F22000270001).}
   \address{ Department of Mathematics and Computer Science, Saint Petersburg State University, Saint Petersburg, 199178, Russia.}
  \email{{\tt pmzlcroak@gmail.com}}

\thanks{The research of P. Mozolyako was supported by the  Ministry of Science and Higher
Education of the Russian Federation, agreement 075-15-2021-602}



\date{}

\maketitle

\pagestyle{myheadings}
\markboth{Two-weight dyadic Hardy inequalities}{Two-weight dyadic Hardy inequalities}

\begin{abstract}
    We present various results concerning the two-weight Hardy's inequality on infinite trees. Our main scope is to survey known characterizations (and proofs) for trace measures, as well as to provide some new ones. Also for some of the known characterizations we provide here new proofs. In particular, we obtain a new characterization based on a new reverse H\"older inequality for trace measures, and one based on the well known Muckenhoupt-Wheeden-Wolff inequality, of which we here give a new probabilistic proof. We provide a new direct proof for the so called isocapacitary characterization and a new simple proof, based on a monotonicity argument, for the so called mass-energy characterization. Furthermore, we introduce a conformally invariant version of the two-weight Hardy's inequality, we characterize the compactness of the Hardy operator, we provide a list of open problems and suggest some possible lines of future research.
\end{abstract}

\newpage

\tableofcontents

\newpage

\section{Introduction}
\subsection{A (very) brief history of Hardy's inequality}
Hardy's inequality \cite{hardy1925} states that

\[ \int_0^\infty \Big( \frac{1}{x}\int_0^x f(y)dy \Big)^p dx \leq (p^*)^p \int_0^\infty f(x)^p dx,\]
for every positive measurable function $f$ and every $p>1$, where $p^*:=p/(p-1)$, and the constant $(p^*)^p$ is optimal. Hardy himself, motivated by the goal of giving a simpler proof of  ``Hilbert's inequality for double series'' \cite{Hardy1920}, was actually primarily interested in the discrete analogue of the above inequality,
\[\label{classicHardy} \tag{Hardy} \quad \sum_{n = 1}^{\infty} \left( \frac{a_1 + \cdots +
   a_n}{n} \right)^p \leq (p^*)^p \sum_{n =
   1}^{\infty} a_n^p, \]
where $a_n$ are positive real numbers. Indeed, for $p=2$, the discrete version was the first one to be proved by Hardy in his earlier paper \cite{hardy1919}. The discrete inequality for general $p$ can either be deduced by the corresponding continuous one or proved directly, as communicated to Hardy by Landau \cite{hardy1925}. Much more information in the fascinating history of the development of Hardy's original inequality can be found in the survey paper \cite{Kufner2006}.

Despite its original purpose it soon became clear that Hardy's inequality and its extensions lie in the heart of the developments in the broad area of harmonic analysis throughout the 20th century, up until modern days. On a macroscopic level that is because Hardy's inequality \emph{is the prototype  of a (weighted) norm inequality for an integration (averaging) operator between $L^p$ spaces.} Averaging operators together with maximal and singular operators are the pillars of harmonic analysis. Under this light a vast number of theorems can be considered as Hardy's type inequalities. Hence, it comes with no surprise that the original Hardy's inequality was later generalized in many different directions.

Tomaselli in \cite{Tomaselli1969} and Talenti in \cite{Talenti1969} made the first steps towards a weighted Hardy's inequality, i.e., an inequality of the form
\begin{equation}\label{Continuous Hardy}
    \int_{0}^\infty \Big( \int_0^x f(y) dy \Big)^p U(x) dx \leq C(p,U,V) \int_0^\infty f(x)^p V(x) dx,
\end{equation}
in which $V,U$ are positive measurable weights. Another way to state the same inequality is to say that the \emph{Hardy operator}, which maps $f$ to its primitive, is bounded from $L^p(\mathbb{R}_+, V(x)dx ) $ to $L^p(\mathbb{R}_+, U(x)dx )$. The first to give a complete characterization of the weights such that the weighted Hardy's inequality holds was Muckenhoupt in \cite{Muckenhoupt1972}. It is worth taking a closer look to Muckenhoupt's condition.

\begin{theorem}[Muckenhoupt]
Let $1\leq p \leq \infty$. There exists a constant $C$ such that \eqref{Continuous Hardy} is true if and only if 
\begin{equation*}
    B:=\sup_{r>0} \Big( \int_r^\infty U(x)^p dx  \Big) \Big( \int_0^r V(x) ^{-p^*}dx \Big)^{p-1}<\infty.
\end{equation*}
Furthermore, if $C$ is the smallest constant such that the inequality holds, then $B\leq C\leq p (p^*)^{p-1} B.$
\end{theorem}

This theorem completes the picture of the weighted Hardy's inequality on $\mathbb{R}_+.$ In the meanwhile several other extensions of Hardy's inequality to different spaces were considered: to higher dimensions, to different metric spaces, to fractional integral operators \cite{Sobolev1938, Stein1958, Muckenhoupt1974}.  Also the weighted problem with different exponents, namely, the boundedness of the Hardy's operator from $L^p(\mathbb{R}_+, V(x)dx ) $ to $L^q(\mathbb{R}_+, U(x)dx )$, with $p$ and $q$ not necessarily coinciding, was extensively studied: we mention Scott \cite{Scott} for the case $1< p\leq q\leq \infty$, Maz'ya \cite{Maz'ya} and Sawyer \cite{Sawyer} for the case $p>q$.

Let us mention that there exists also a different stream of research which is dedicated to finding weights which are optimal, in an appropriate sense, for weighted Hardy's inequalities, both in the continuous setting (\cite{Devyver2014, Berchio2020}) and in the setting of graphs (\cite{pinchover, berchio}).


\subsection{The two-weight Hardy's inequality on trees}
In the present paper we focus on a generalized version of the discrete Hardy's inequality \eqref{classicHardy}, the two-weight Hardy's inequality on trees. A particular case of this inequality is the so called two-weight dyadic Hardy's inequality, presented in Section \ref{sec:dyadic} in the Appendix. In that section we will also make clear the intuitive fact that the classical inequality \eqref{classicHardy} is a special case of the inequality on trees, and we will provide some examples of application of the dyadic inequality in Complex Analysis. In the paper, however, we are able to work in the generality of the two-weight Hardy's inequality on trees presented in this section. Before stating such  an inequality we need to introduce some pieces of notation.

A \it tree \rm $T=(V,E)$ is a simple connected graph with no cycles, where $V$ is a finite or countable set of vertices and
$E$ is the set of edges. The number of edges sharing each vertex is assumed to be finite, but we don't require further restrictions. We will consistently use Greek letters for edges and Latin letters 
for vertices and boundary points, to be defined below. Henceforth we will identify the tree $T$ with its vertices $V$. We fix arbitrarily a \it root \rm vertex $o$ and we assume there exists a \it pre-root \rm vertex $o^*$ which is connected to $o$ and to no other vertex. We denote by $\omega$ the edge connecting $o$ and $o^*$ and we call it the \it root edge \rm.

A fundamental property of trees is that for any couple of vertices $x,y\in T$ there exists a unique geodesic connecting the two, that is, a unique minimal sequence of pairwise connected vertices containing $x$ and $y$. We write $[x,y]$, or equivalently $[y,x]$, for the (unique) set of edges connecting pairs of points in the geodesic joining $x$ and $y$. The \it confluent \rm of $x$ and $y$ is the vertex $x\wedge y$ such that $[o^*,x\wedge y]=[o^*,x]\cap[o^*,y]$. The edge-counting distance on $T$ is given by $d(x,y)=\sharp[x,y]$. We use the same symbol for the vertex-counting distance on $E$, defined in the obvious way, and abbreviate $d(\alpha)$ for $d(\alpha,\omega)$. If we assign to each edge a weight $\sigma(\alpha)>0$, we can define the associated distance $d_\sigma$ by adding the weights 
on the edges of a paths, instead of counting edges. Due to the elementary topology of a tree,
the new metric has the same geodesics.

 If $\alpha$ is an edge, we denote by $b(\alpha)$ its endpoint vertex which is closest to $o^*$ and by $e(\alpha)$ the furthest. Observe that for any vertex $x$ there exists a unique edge $\alpha$ with $e(\alpha)=x$, while there are possibly many edges with $b(\alpha)=x$. For each vertex $x$ we define the \it predecessor \rm of $x$ as its unique neighbor vertex $p(x)$ which is closer than $x$ to $o^*$, and we denote by $s(x)$ the set of the remaining neighbors, the \it children \rm on $x$.
 Predecessors and children may be defined also for edges in the very same way.

The choice of a root induces a partial order on $T$ and $E$. 
For edges, we write $\alpha\subseteq\beta$ if $\beta\in[o^*, e(\alpha)]$, and for vertices we write $x\subseteq y$ if
$[o^*,y]\subseteq[o^*,x]$. We will also write $x\subseteq\alpha\subseteq y$, comparing vertices and edges, with the obvious
meaning. The notation $\subseteq$ is a reminder that trees often come from dyadic decompositions of metric spaces (see Section \ref{sec:dyadic}).

The \it boundary \rm of the rooted tree
$(T,\omega)$, denoted by $\partial T$, is the set of the maximal geodesics emanating from $o^*$ \footnote{If a maximal geodesic is of finite lenght, i.e., it ends in a \textit{leaf} vertex $x$, we identify such a geodesic with $x$, which is then a boundary point for $T$.}.
In the infinite tree case, we always assume that all maximal geodesics 
starting at $o^*$ are infinite, so that $\partial T\cap T=\emptyset$. We set $\overline{T}:=T\cup \partial   T$.

The most important geometric objects we deal with are the \it successor sets \rm $S(\alpha)$ of an edge $\alpha$,
$S(\alpha)=\{x\in \overline{T}:\ [o^*,x]\ni\alpha\}$. We also write $S(x)$ for $S(\alpha)$ when $x=e(\alpha)$. We remark that $\overline{T}$ is compact with respect to the topology generated by the family $\lbrace S(\alpha)\rbrace_{\alpha\in E}$ of successor sets and singletons $\{x\}_{x\in T}$, from which $\overline{T}$ is often referred to as the \it standard compactification \rm of $T$.

The following picture should help the reader to visualize and summarize some of the fundamental notation just introduced.

\begin{figure}[H]

\scalebox{0.75}{%
\tikzstyle{every node}=[font=\large]

\begin{tikzpicture}

\node[shape=circle,draw=black, fill=black, scale=0.5, label=left:{$o^*$}, label=right:{ \ pre-root}](O) at (0,2) {};

\node[shape=circle,draw=black, fill=black, scale=0.5, label=left:{$o$ \ }, label=right:{ \ root}](A) at (0,0) {};

\node[shape=circle,draw=black, fill=black, scale=0.5](Aa) at (-5,-2) {};
\node[shape=circle,draw=black, fill=black, scale=0.5, label={left:{$b(\alpha)$}}, label={right:{$p(x)$}}](Ab) at (0,-2) {};
\node[shape=circle,draw=black, fill=black, scale=0.5](Ac) at (5,-2) {};

\node[shape=circle,draw=black, fill=black, scale=0.5](Aaa) at (-7,-4) {};
\node[shape=circle,draw=black, fill=black, scale=0.5](Aab) at (-5,-4) {};
\node[shape=circle,draw=black, fill=black, scale=0.5](Aac) at (-3,-4) {};

\node[shape=circle,draw=black, fill=black, scale=0.5, label={left:{$e(\alpha)$}}, label={right:{$x$}}](Aba) at (-1,-4) {};
\node[shape=circle,draw=black, fill=black, scale=0.5](Abb) at (1,-4) {};

\node[shape=circle,draw=black, fill=black, scale=0.5](Aca) at (4,-4) {};

\node[shape=circle,draw=black, fill=black, scale=0.5](Acb) at (6,-4) {};

\node[shape=circle,draw=black, fill=black, scale=0.5](Aaaa) at (-7.5,-6) {};
\node[shape=circle,draw=black, fill=black, scale=0.5](Aaab) at (-6.5,-6) {};

\node[shape=circle,draw=black, fill=black, scale=0.5](Aaba) at (-5.5,-6) {};
\node[shape=circle,draw=black, fill=black, scale=0.5](Aabb) at (-4.5,-6) {};

\node[shape=circle,draw=black, fill=black, scale=0.5](Aaca) at (-3.5,-6) {};
\node[shape=circle,draw=black, fill=black, scale=0.5](Aacb) at (-2.5,-6) {};

\node[shape=circle,draw=black, fill=black, scale=0.5](Abaa) at (-1.5,-6) {};
\node[shape=circle,draw=black, fill=black, scale=0.5](Abab) at (-0.5,-6) {};

\node[shape=circle,draw=black, fill=black, scale=0.5](Abba) at (0.5,-6) {};
\node[shape=circle,draw=black, fill=black, scale=0.5](Abbb) at (1.5,-6) {};

\node[shape=circle,draw=black, fill=black, scale=0.5](Acaa) at (3.5,-6) {};
\node[shape=circle,draw=black, fill=black, scale=0.5](Acab) at (4.5,-6) {};

\node[shape=circle,draw=black, fill=black, scale=0.5](Acba) at (5.5,-6) {};
\node[shape=circle,draw=black, fill=black, scale=0.5](Acbb) at (6.5,-6) {};

\node[draw=none, fill=none](Aaaaa) at (-7.5,-7) {};
\node[draw=none, fill=none](Aaaba) at (-6.5,-7) {};

\node[draw=none, fill=none](Aabaa) at (-5.5,-7) {};
\node[draw=none, fill=none](Aabba) at (-4.5,-7) {};

\node[draw=none, fill=none](Aacaa) at (-3.5,-7) {};
\node[draw=none, fill=none](Aacba) at (-2.5,-7) {};

\node[draw=none, fill=none](Abaaa) at (-1.5,-7) {}; 
\node[draw=none, fill=none](Ababa) at (-0.5,-7) {};

\node[draw=none, fill=none](Abbaa) at (0.5,-7) {};
\node[draw=none, fill=none](Abbba) at (1.5,-7) {};

\node[draw=none, fill=none](Acaaa) at (3.5,-7) {};
\node[draw=none, fill=none](Acaba) at (4.5,-7) {};

\node[draw=none, fill=none](Acbaa) at (5.5,-7) {};
\node[draw=none, fill=none](Acbba) at (6.5,-7) {};

  
  \path [- ] (O) edge node[left]{$\omega$} node[right]{ \ root edge}(A); 
  \path [- ] (A) edge node[left] {} (Aa); 
  \path [-] (A) edge node[left] {} (Ab); 
  \path [-] (A) edge node[right] {} (Ac); 

 \path [- ] (Aa) edge node[left] {} (Aaa); 
 \path [-] (Aa) edge node[left] {} (Aab); 
 \path [-] (Aa) edge node[left] {} (Aac); 

\path [-,black ] (Ab) edge node[left] {$\alpha$} (Aba);
 \path [-] (Ab) edge node[left] {} (Abb); 
 
 \path [- ] (Ac) edge node[left] {} (Aca); 
 \path [-] (Ac) edge node[right] {} (Acb); 
 
 \path [- ] (Aaa) edge node[left] {} (Aaaa);  \path [- ] (Aaa) edge node[left] {} (Aaab);
 \path [- ] (Aab) edge node[left] {} (Aaba);  \path [- ] (Aab) edge node[left] {} (Aabb);
 \path [- ] (Aac) edge node[left] {} (Aaca);  \path [- ] (Aac) edge node[left] {} (Aacb);
\path [-,black ] (Aba) edge node[left] {} (Abaa);  \path [-,black ] (Aba) edge node[left] {} (Abab);
 \path [- ] (Abb) edge node[left] {} (Abba);  \path [- ] (Abb) edge node[left] {} (Abbb);
 \path [- ] (Aca) edge node[left] {} (Acaa);  \path [- ] (Aca) edge node[left] {} (Acab);
 \path [- ] (Acb) edge node[left] {} (Acba);  \path [- ] (Acb) edge node[right] {} (Acbb);
 
 \path [dashed] (Aaaa) edge node[left] {} (Aaaaa);  \path [dashed] (Aaab) edge node[left] {} (Aaaba);
 \path [dashed] (Aaba) edge node[left] {} (Aabaa);  \path [dashed] (Aabb) edge node[left] {} (Aabba);
 \path [dashed] (Aaca) edge node[left] {} (Aacaa);  \path [dashed] (Aacb) edge node[left] {} (Aacba);
 \path [dashed] (Abaa) edge node[left] {} (Abaaa);  \path [dashed] (Abab) edge node[left] {} (Ababa);
 \path [dashed] (Abba) edge node[left] {} (Abbaa);  \path [dashed] (Abbb) edge node[left] {} (Abbba);
 \path [dashed] (Acaa) edge node[left] {} (Acaaa);  \path [dashed](Acab) edge node[left] {} (Acaba);
 \path [dashed] (Acba) edge node[left] {} (Acbaa);  \path [dashed] (Acbb) edge node[left] {} (Acbba);
 \path [dashed] (Acaa) edge node[left] {} (Acaaa);  \path [dashed] (Acab) edge node[left] {} (Acaba);
 
 \filldraw [fill=blue,draw opacity=0,,fill opacity=0.5] (-1,-4) -- (-2,-8) -- (0,-8) -- cycle;
 \node[draw=none, fill=none] at (1.5,-7.5) {\color{blue} $S(\alpha)=S(x)$}; 

 \draw[dashed] (-8,-8) -- (6.8,-8);
 \node[draw=none, fill=none] at (7.2,-8) {$\partial T$}; 
 
  \node[draw=none, fill=none] at (-1,-8.5) {\color{blue} $\partial S(\alpha)$};

\end{tikzpicture}
}
\caption{The infinite rooted tree.}
    \label{fig:1}
\end{figure}
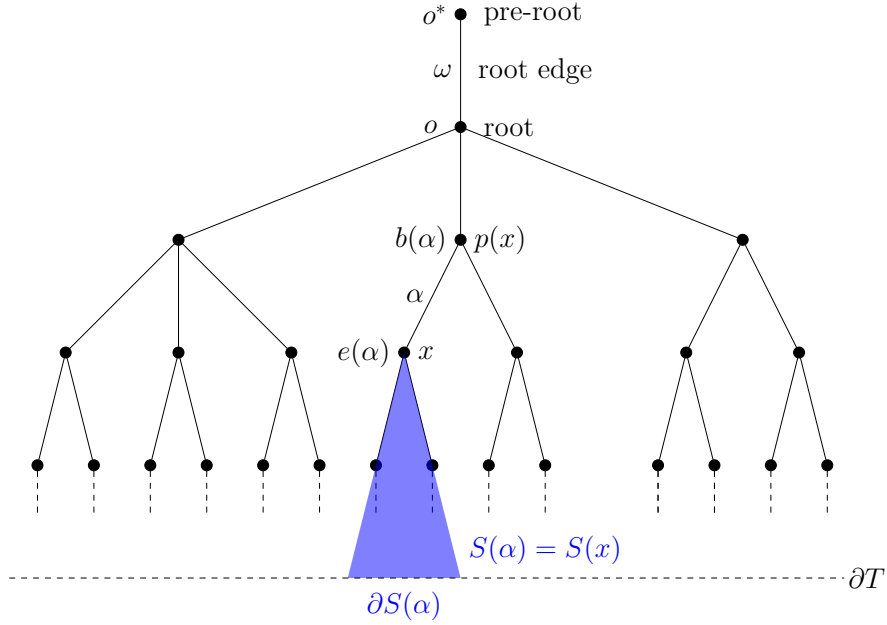

We are now ready to introduce our main object of study, the two-weight Hardy's inequality on trees. Define the \textit{Hardy operator} $\III$ as the operator mapping a function $\varphi:E\to \RR_+$ to
\begin{equation}\label{HardyOperator}
\III\varphi(x)=\sum_{\alpha\in[o^*,x]}\varphi(\alpha), \quad x\in \overline{T},
\end{equation}
provided that the sum converges. Let $\pi:E\to \RR_+$ be some fixed edge weight, and $\mu\ge0$ a Borel measure on $\overline{T}$ and write $\ell^p(\pi)=\ell^p(E,\pi)$ and $L^p(\mu)=L^p(\overline{T},\mu)$. The main problem which is discussed 
here is whether the following \textit{two-weight Hardy's inequality on $T$} holds,
\begin{equation}\label{HardyInequality}\tag{H}
\int_{\overline{T}}\III\varphi(x)^pd\mu(x)\le [\mu] \sum_{\alpha\in E}\varphi(\alpha)^p\pi(\alpha),
\end{equation}
i.e., if $\III : \ell^p(\pi)\to L^p(\mu)$  is bounded, and what can be said about the constant $[\mu]=\|\III\|_{\ell^p(\pi)\to L^p(\mu)}^p$.

Any $\mu$ satisfying \eqref{HardyInequality} (for some fixed $\pi$, $p$) with finite constant $[\mu] $ will be called a \textit{trace measure}, and $[\mu]$ is the \it Carleson measure norm, \rm or \it trace norm, \rm of $\mu$.
We always consider $\pi$ as a fixed, geometric object so that $[\mu]$ depends on $\pi$ and $p$ as well. Also, the quantity $[\cdot]$ is chosen so to be sublinear functional of $\mu$:
\begin{equation*}
   [t\mu]=t[\mu] \text{ if }t\ge0,\text{ and } [\mu+\nu]\le[\mu]+[\nu].
\end{equation*}

The inequality \eqref{HardyInequality} has a natural corresponding dual inequality, which we will be of fundamental use in the sequel. For $\psi:V\to\mathbb{R}$, we set
\begin{equation}\label{DualHardyOperator}
 \III^*_\mu\psi(\alpha)=\III^*(\psi d\mu)(\alpha):=\int_{S(\alpha)}\psi d\mu \quad \alpha\in E,
\end{equation}
By obvious duality,
\begin{equation}\label{HardyDual}
\langle\III\varphi,\psi\rangle_{L^2(\mu)} = \langle\varphi,\III^*(\psi d\mu)\rangle_{\ell^2}, 
\end{equation}
Thus, \eqref{HardyInequality} is equivalent to
\begin{equation}\tag{$\text{H}^*$} \label{HardyInequalityDual}
 \sum_{\alpha}\pi(\alpha)^{1-\pp}(\III^*_\mu\psi(\alpha))^\pp\le [\mu]^{\pp-1}\int_{\overline{T}}\psi^\pp d\mu,
\end{equation}
which will refer to as the \emph{dual Hardy's inequality.}

We reserve a specific symbol for a particular edge weight which will show up often in the paper; we write $|\alpha|$ for the unique edge weight satisfying the recursive formula $|\alpha|=|p(\alpha)|/q(e(\alpha))$, normalized such that $|\omega|=1$. On the homogeneous tree where each vertex has $q+1$ neighbors, $|\alpha|=q^{-d(\alpha)}$. Observe that if the tree arises from a dyadic decomposition of the unit interval (see Section \ref{sec:dyadic}) then $|\alpha|$ is the Lebesgue measure of the corresponding subinterval in the decomposition. In full generality we will then call Lebesgue measure on the boundary of $T$ the measure $dx$ defined by
\begin{equation}
    \int_{\partial S(\alpha)}dx=|\alpha|, \quad \alpha \in E.
\end{equation}

Characterizations of the triples $(p, \pi, \mu)$ for which \eqref{HardyInequality} holds have
been know since more than twenty years (\cite{Sawyer, Arcozzi2002, AHMV, Cavina2020}). The main
scope of this expository paper is surveying old characterizations together with their proofs, provide some new proofs and also present some completely new characterizations. We point out that even the known proofs are
here adapted to work in our general framework, since they have all been originally
given for the dyadic homogeneous tree only.

Our motivations are manifold. First, it is interesting and instructive to see
the diverse machinery which can be employed in the solution of the problem.
Second, we will see that the many conditions characterizing $\mu' s$ for which
the Hardy's inequality holds are rather different one from the other, and their
equivalence is a collection of interesting mathematical facts by itself.
Third, as research moves to uncharted territories, such as multi-parameter
dyadic Hardy inequalities, it is useful to have a place were different
techniques are surveyed in a unified framework. We will mention some more
recent areas of investigation, providing some results and mentioning a number
of open problems.

\subsection{Description of contents}

We proceed now to a detailed description of the main theorems of this paper. The main goal of Sections \ref{sec:iso} and \ref{sec:me} is to prove the following theorem. 

\begin{theorem}\label{main theorem} Let $T$ be a locally finite connected tree, $1<p<\infty$ and $\pi$ a non negative weight on the edges. Then for a positive Borel measure $\mu$ on $\overline{T}$ the following are equivalent. 

\begin{itemize}
    \item[(i)] The two-weight Hardy's inequality \eqref{HardyInequality} holds with best constant $[\mu]<  \infty.$
    \vspace{0.2cm}
    \item[(ii)] The {\it mass-energy condition} holds; 
    \begin{equation}\label{EqMassEnergy} \tag{ME}
      \sup_{\alpha \in E } \mu(S(\alpha)) ^{-1}\sum_{\beta\subseteq\alpha}\pi(\beta)^{1-\pp}\mu(S(\beta))^\pp := [[\mu]]^{\pp-1}<\infty.
      \end{equation}
      \item[(iii)] The isocapacitary condition holds;
       \begin{equation}\label{eq:isocapacitary}\tag{ISO}
         \sup_{\alpha_1, \dots, \alpha_n \in E}  \Capac_{p,\pi}\Big( \bigcup_{i=1}^n S(\alpha_i) \Big) ^{-1}   \sum_{i=1}^n\mu(S(\alpha_i)) := [[\mu]]_c , < \infty,
     \end{equation}
     where the capacity $\Capac_{p,\pi}$ is the one defined in Section \ref{subsec:Some potential theory}.
\end{itemize}
Furthermore, the following inequalities hold,
\[ [[\mu]]_c \leq [\mu] \leq 2^p [[\mu]]_c, \quad [[\mu]] \leq [\mu] \leq p^{p^*} [[\mu]]. \]
\end{theorem}


As we will show later, it is easy to  prove that \eqref{EqMassEnergy} and \eqref{eq:isocapacitary} imply \eqref{HardyInequality}; the main issue in Theorem \ref{main theorem} is to prove the reverse implications. Section \ref{sec:iso} is dedicated to developing a potential theory on the tree and proving the isocapacitary characterization , i.e., the equivalence of (i) and (iii) in Theorem \ref{main theorem}. This equivalence was first proved in \cite{ARS2008}, though in an indirect way, passing through the mass energy condition. Here we give a new direct proof that  (iii) implies (i), which builds on ideas developed by Maz'ya in the continuous setting. The main tool is a \emph{strong capacitary
inequality} \cite{Mazya1973, Adams1976}, which in tree language takes the form,
  \[\tag{Cap} \label{strongCap} \sum_{k = - \infty}^{+ \infty} 2^{p k} \Capac_{p,\pi} (\{x : \III \varphi (x) > 2^k\})
     \leq \frac{2^p}{2^p-1} \| \varphi  \|_{\ell^p (E,\pi)}^p , \]
and has an elementary proof.

In Section \ref{sec:me} we turn to the \emph{mass--energy condition} introduced in \cite{Arcozzi2002}.  We will give three proofs of its equivalence  with  \eqref{HardyInequality}: one based on
maximal functions (originally proved for the homogeneous tree in \cite{ARS2007}), one relying on a
simple monotonicity argument, which is new and the simplest available at the
moment, but only works for $p = 2$, and a very recent one using a
Bellman function argument (originally proved for the homogeneous tree in \cite{AHMV, Cavina2020}).

The advantage of the isocapacitary condition \eqref{eq:isocapacitary}  over the mass--energy
condition \eqref{EqMassEnergy} is that the measure $\mu$ appears on the left
hand side only. It is obvious from it, for instance, that if $\nu \leq
\mu$, then $[[\nu]]_c\leq [[\mu]]_c$. On the other hand, the mass--energy
condition only has to be verified on single intervals, and not on arbitrary
unions of them.

In some simple and particularly common cases of weights and trees, some further characterizations for trace measures can be proved to hold, lengthening the list of equivalent conditions summed up in Theorem \ref{main theorem}. Sections \ref{sec:holder} and \ref{sec:MW} are dedicated to two different such additional characterizations.

More precisely, in Section \ref{sec:holder}, we restrict our attention to the case $\pi = 1 $ and $ p=2 $, and prove that in this case the Hardy's inequality is equivalent to a one parameter family of conditions. We prove the following theorem.
\begin{theorem}
The Hardy's inequality
\begin{equation*}
\int_{\overline{T}}\III\varphi(x)^2d\mu(x)\le [\mu] \sum_{\alpha\in E}\varphi(\alpha)^2,
\end{equation*}
is equivalent to
 \[\tag{s-Testing} \sup_{\alpha \in E} \dashint_{S(\alpha)} \Big( \int_{S(\alpha)} d(x\wedge y) d\mu(x) \Big)^s d\mu(y): = [[\mu]]_s < \infty.\]
for some (equivalently, for every) $s\geq1$.
\end{theorem}
It can be readily verified that the $s-$testing condition is stronger than the mass--energy condition. The aim of the section is to prove that in fact they are all equivalent to the mass--energy condition which, in a sense, amounts to say that Carleson measures satisfy a reverse H\"older inequality, see Theorem \ref{CalderonZygmound}. The results in this section are new and the techniques employed are partially inspired by the work of Tchoundja \cite{Tchoundja2008}.

In Section \ref{sec:MW} we provide another characterization of trace measures which holds (for any $p$) for a family of edge weights (depending on $p$) on homogeneous trees, and it can be generalized to trees having Ahlfors regular boundary \cite[Section 3]{Arcozzi2013}. More precisely we prove the following.
\begin{theorem}
  Let $T$ be a homogeneous tree, $1<p<\infty$, $0<s<1$ and $\pi(\alpha) = |\alpha|^{\frac{1-p^*s}{1-p^*}}$. Then, the following conditions are both equivalent to the Hardy's inequality \eqref{HardyInequality}
\begin{itemize}
    \item[(i)] $\displaystyle
        \sup_{\alpha\in E}  \mu(S(\alpha)) ^{-1}\int_{\partial T_\alpha} \left( \sum_{\beta \supset x} \frac{\mu (S
     (\beta))}{|\beta|^s} \right)^{p^*} d x<\infty$;
 \vspace{0.2cm}
    \item[(ii)] 
    $\displaystyle
      \sup_{\alpha\in E}\mu(S(\alpha)) ^{-1}\int_{\partial T} \left( \sup_{\beta \supset x}
     \frac{\mu (S (\beta))}{|\beta|^s} \right)^{p^*} d x<\infty$.
\end{itemize}
\end{theorem}

Such a characterization is an immediate consequence of the
Muckenhoupt--Wheeden inequality, which in this case reads as follows: for any $1<p<\infty$, $0<s<1 $  and for any measure $\mu$ on a homogeneous tree $T$,

 \begin{equation*}
      \int_{\partial T} \left( \sum_{\alpha \supset x} \frac{\mu (S
     (\alpha))}{|\alpha|^s} \right)^{p^*} d x \approx \sum_{\alpha } \frac{\mu (S(\alpha))^{p^*}}{|\alpha|^{p^*s-1}} \approx  \int_{\partial T} \left( \sup_{\alpha \supset x}
     \frac{\mu (S (\alpha))}{|\alpha|^s} \right)^{p^*} d x. 
  \end{equation*}
Indeed, it is easy to recognize in the middle term the form that the sum appearing in the mass--energy condition gets for the particular choice of $\pi(\alpha) = |\alpha|^{\frac{1-p^*s}{1-p^*}}$. Hence, for this family of weights on homogeneous trees we have an alternative characterization of Carleson measures.  The above choice of $\pi$ is of particular interest because of the connection with the theory of Bessel potentials in $\mathbb{R}^n$, see Section \ref{subsec:Bessel} in Appendix.

Besides providing a different characterization of Carleson measures, the story of this inequality is itself interesting. That the term one the left is comparable to that on the right was proved in a non--dyadic language
 by Muckenhoupt and Wheeden in  \cite[Theorem 1]{Muckenhoupt1974}. They attribute the idea of the proof, which is a textbook example of a good lambda inequality, to Coifman and Fefferman \cite{Coifman1974}.
Unaware of it, Wolff gave a
wholly different proof \cite[Theorem 1]{Hedberg1983} that the term in the center
is comparable to the term on the left. Independently of this, the first author,
 Rochberg and Sawyer gave a different proof that the central term is bounded by
the one on the right \cite{Arcozzi2002}. In this section we will give a new proof of the ``Wolff inequality'' based
on a probabilistic argument. The new proof has the advantage that it extends more easily to different settings. We will return on that in future work.

In Section \ref{sec: conf inv} we discuss 
a conformally invariant version of \eqref{HardyInequality} on the homogeneous tree. The left hand side in \eqref{HardyInequality} evidently depends on the arbitrary choice of a root in our tree. To remedy this we propose the following  modified inequality 
\begin{equation} \tag{CH}
    \int_{\overline{T}} \Big| \mathcal{I} f (x) - \frac{1}{\mu(\overline{T})} \int_{\overline{T}} f d\mu \Big|^2 d \mu (x) \leq [\mu]_{inv}
   \sum_\alpha f (\alpha)^2.
\end{equation}
We prove that it is ``invariant'', i.e., if $\Psi$ is a tree automorphism, then $[\Psi_*\mu]_{inv}=[\mu]_{inv}$, and that it is surprisingly equivalent to \eqref{HardyInequality}. We also provide a sharp estimate of the quantity $[\mu]_{inv}$ in terms of capacity (Theorem \ref{thm:hardy equiv hardy inv}). All results in this section appear here for the first time.

In Section \ref{sec:miscellaneous} we collect some miscellaneous results on the topic. They are all new. First we prove that a ``vanishing'' version of the mass--energy condition characterizes the compactness of the Hardy operator (Theorem \ref{comppp}). With a similar reasoning one can obtain an equivalent isocapacitary type vanishing condition (Theorem \ref{comp-cap}). We then discuss another very natural and easily determined necessary condition for \eqref{HardyInequality} to hold, the \it simple box-type condition\rm,
\begin{equation}\label{CarlesonBox} \tag{SB}
 \sup_{\alpha \in E}\mu(S(\alpha)) \left(\sum_{\beta\supseteq\alpha}\pi(\beta)^{1-\pp}\right)^{p-1} :=[[\mu]]_{sc}<+\infty.
\end{equation}
In contrast to the mass-energy and the isocapacitary conditions, however, for the many relevant weights and trees \eqref{CarlesonBox} is not sufficient, see Example \ref{counterexample}. We end the section by providing two easy examples of trees where the potential theory degenerates in two opposite ways.

In Section \ref{sec:variation} we enlarge our horizon, including some dyadic structures
variously related to the Hardy operator, or to its applications. Most of this
territory is uncharted, only a few results are known, investigation is still
in its infancy, and there is a high potential for applications to harmonic
analysis, holomorphic function theory, and more. This section is essentially a description of the few results that are known in the literature. We omit most of the proofs. In Section \ref{sec:rk} we introduce the viewpoint of
reproducing kernel Hilbert spaces, which provides a unified view of the
preceding inequalities and is instrumental to state the problem of Hardy-type
inequalities for quotient structures in Section \ref{sec:quotient}.  In Section \ref{sec:product} we briefly account on the topic of
Hardy inequalities on poly-trees. This is a new area of research where very little is known. Recently it has attracted a lot of interest because of the applications to function theory in the poly--disc.

We end the paper with an Appendix where we include the discussion on the model case of the purely dyadic Hardy's inequality (Section \ref{sec:dyadic}) and a comparison of the potential theory we use in the paper with that arising from Bessel's potentials (Section \ref{subsec:Bessel}.

In the text we mention a number of open problems that we think are interesting and deserve further attention.

\section{Potential theory on trees and the isocapacitary characterization}\label{sec:iso}

This section is dedicated to prove the equivalence of (i) and (iii) in Theorem \ref{main theorem}. Section \ref{subsec:Some potential theory} introduces the potential theory which we need to define a $p$-capacity on the tree, while the actual proof is given in the subsection \ref{sec:SCI}.

\subsection{Potential theory}\label{subsec:Some potential theory}

We define a potential theory following Adams and Hedberg's axiomatic approach \cite{Adams_book}. Other approaches are also possible, see for example \cite{Soardi1994}.
Consider the compact Hausdorff space $\overline{T}$, and make $E$ into a measure space by endowing it with the  measure associated to a weight $\sigma:E\to\mathbb{R}_+$. We introduce the kernel $k:\overline{T}\times E\to \mathbb{R}_+$, given by the characteristic function $k(x,\alpha)=\chi_{\lbrace\alpha \supset x \rbrace}(x,\alpha)$. Observe that $k(\cdot,\alpha)$ is continuous on $\partial T$, since $\partial S(\alpha)$ is open.

Given a function $ \varphi: E \rightarrow \mathbb{R}_+ $, we define the \textit{potential} of $\varphi$, $\III_\sigma\varphi: \overline{T} \rightarrow \mathbb{R}_+\cup\{+ \infty \} $, by
\begin{equation*}
  \III_\sigma\varphi(x)=\sum_{\alpha} k(x,\alpha)\varphi(\alpha)\sigma(\alpha)=\sum_{E\ni\alpha \supset x }\varphi(\alpha)\sigma(\alpha).  
\end{equation*}

The \textit{co-potential} of a function $\psi:\overline{T}\to\mathbb{R}_+$ with respect to a positive Borel measure $\mu$ on $\overline{T}$ is defined as the edge function
\begin{equation*}
    \III^\ast_\mu\psi(\alpha)=\int_{\overline{T}}k(x,\alpha)\psi(x)d\mu(x)=\int_{S(\alpha)}\psi(x)d\mu(x), \quad \alpha\in E.
\end{equation*}

The co-potential of $\mu$ is intended to be $\III^\ast_\mu(\alpha):=\III^\ast_\mu 1(\alpha)=\mu(S(\alpha))$. Observe that if $\psi\in L^1(\overline{T},\mu)$, by Fubini's theorem we have $\langle \III_\sigma\varphi,\psi\rangle_{L^2(\overline{T},\mu)}=\langle \varphi, \III^\ast_\mu\psi\rangle_{\ell^2(E,\sigma)}$.

For a H\"older dual pair of exponents $p, \pp $  we can further associate to the measure $\mu$ a \textit{nonlinear Wolff's potential}, $V_p^{\mu,\sigma}(x)=\III_\sigma(\III^\ast_\mu)^{p^*-1} (x)$, $x\in \overline{T}$. More explicitly, one has
\begin{equation*}
    V_p^{\mu,\sigma}(x)=\sum_{\alpha}k(\alpha,x)\sigma(\alpha)\Big(\int_{\overline{T}}k(\alpha,y)d\mu(y)\Big)^{p^*-1}=\sum_{\alpha\in [o^*,x]}\sigma(\alpha)\mu(S(\alpha))^{p^*-1}.
\end{equation*}
In the linear case $p=2$ the sum and the integral can be switched and the above potential expressed as
\begin{equation*}
    V_2^{\mu,\sigma}(x)=\int_{\overline{T}}\sum_{\alpha}k(\alpha,x\wedge y)\sigma(\alpha)d\mu(y)=\int_{\overline{T}}\sum_{\alpha\supset x\wedge y }\sigma(\alpha)d\mu(y).
\end{equation*}
The \textit{$p-$energy} of the charge distribution $\mu$ is given by
\begin{equation*}
    \mathcal{E}^\sigma_p(\mu)=\int_{\overline{T}}V_p^{\mu,\sigma}(x)d\mu(x).
\end{equation*}
If the energy is finite, by Fubini's theorem it holds
\begin{equation*}
    \mathcal{E}^\sigma_p(\mu)=\Vert \III^\ast_\mu\Vert_{\ell^{p^*}(E,\sigma)}^{p^*}=\sum_{\alpha}\mu(S(\alpha))^{p^*}\sigma(\alpha).
\end{equation*}
We define the capacity of a closed subset $A\subseteq\overline{T}$ as
\begin{eqnarray*}
  \tmop{Cap}_p^\sigma (A) & = & \inf \left\{ \Vert \varphi\Vert_{\ell^p(E,\sigma)}^p :
  \ \varphi \geq 0, \ \III_\sigma \varphi (x)\geq 1
  \text{ on } A \right\}\\
  & = &  \sup \left\{ \frac{\mu (A)^p}{\mathcal{E}^\sigma_p (\mu)^{p -
   1}} : \ \text{supp} (\mu) \subseteq A \right\}.
\end{eqnarray*}
The equality between the first and second line above is given by a classical theorem in potential theory that can be found, for instance, in \cite{Adams_book}. The $(p,\sigma)-$equilibrium function for $A$ is the unique function $\varphi$ satisfying $\III_\sigma \varphi = 1$, $\tmop{Cap}_p^\sigma-$quasi everywhere on $A$, and $\tmop{Cap}_p^\sigma (A)=\Vert \varphi\Vert_{\ell^p(E,\sigma)}^p$. Similarly, one defines the $(p,\sigma)-$equilibrium measure for $A$ as the unique measure probability measure such that $\tmop{Cap}_p^\sigma (A)=\mu(A)$. Two of the authors recently found a characterization for equilibrium measures on trees \cite{equilibrium}.

Observe that the boundedness of $\III_\sigma: \ell^p(\sigma)\to L^p(\mu)$ is equivalent to that of the Hardy operator $\III: \ell^p(\pi)\to L^p(\mu)$, under the correspondence $\pi(\alpha)=\sigma(\alpha)^{1-p}$. Following this paradigm, we can translate the natural potential theoretic objects in the $\pi-$dictonary, which turns out to be more adjusted to our scopes:
\begin{eqnarray*}
    & &\mathcal{E}^\sigma_p(\mu)=\mathcal{E}_{p,\pi}(\mu):=\sum_{\alpha}\mu(S(\alpha))^{p^*}\pi(\alpha)^{1-p^*},\\
    & &\tmop{Cap}_p^\sigma (A)=\tmop{Cap}_{p,\pi} (A)=\inf \left\{ \Vert \varphi\Vert_{\ell^p(E,\pi)}^p :
  \ \varphi \geq 0, \ \III \varphi (x)\geq 1
  \text{ on } A \right\}.
\end{eqnarray*}

\subsection{Strong Capacitary Inequality and the isocapacitary characterization}\label{sec:SCI}

We are ready to prove the isocapacitary characterization for trace measures, that is, the equivalence of (i) and (iii) in Theorem \ref{main theorem}. As it is common, we will prove it passing through the so called capacitary strong inequality. There are various versions of such inequality. Here we naturally treat the case of a tree, a proof for a large class of kernels in the continuous case can be found in \cite[Theorem 7.1.1]{Adams_book}.

\begin{theorem}
  [Capacitary Strong Inequality]
  Let $1<p<\infty$ and $\varphi : E \to \mathbb{R}_{+}$
  \[ \sum_{k = - \infty}^{+ \infty} 2^{p k} \Capac_{p,\pi} (\{x : \III \varphi (x) > 2^k\})
     \leq \frac{2^p}{2^p-1} \| \varphi  \|_{\ell^p (E,\pi)}^p . \]
\end{theorem}

\begin{proof}
Set $\Omega_k = \{ x \in T : \III \varphi (x) > 2^k \}$, $\partial \Omega_k = \{ x
\in \Omega_k : \III \varphi (p(x)) \leq 2^k \}$. For given $k$, let $\varphi_k (\alpha) =2^{k-m(\alpha)}\varphi (\alpha)$ if $e(\alpha)\in \partial\Omega_k\cap\partial\Omega_{k+1}$, where $m(\alpha)=m$ is the largest integer for which $e(\alpha) \in \partial\Omega_{m + 1}$, and $\varphi_k (\alpha)=\varphi (\alpha)\cdot \chi_{\Omega_k \setminus \Omega_{k + 1}} (e(\alpha))$ otherwise. Let $x\in \partial\Omega_{k+1}$ and $y=[o,x]\cap\partial\Omega_k$. Then, if $x\neq y$, we have
\[ \III \varphi_k (x) = \III \varphi (x) - \III \varphi (p(y)) > 2^{k + 1} - 2^k = 2^k , \]
while if $x=y$ and $\alpha$ is the unique edge such that $x=e(\alpha)$, we have $\varphi (\alpha)=\III \varphi (x) - \III \varphi (p(x))>2^{m + 1} - 2^m$ and therefore
\[ \III \varphi_k (x) = \varphi_k (\alpha)=2^{k-m}\varphi (\alpha)>2^k.\]
Hence, $2^{ - k} \varphi_k$ is a testing function for $\text{Cap}_{p,\pi} (x : \III \varphi
(x) > 2^k)$. Summing and using that using that the supports of the $\varphi_k' s$ can meet only
at points belonging to multiple boundaries,
\begin{eqnarray*}
  \sum_{k\in\mathbb{Z}} 2^{p k} \text{Cap}_{p,\pi} (\Omega_{k}) & \leq & \sum_{k\in\mathbb{Z}} 2^{p k}
  \sum_\alpha 2^{- p k} | \varphi_k (\alpha) |^p \pi (\alpha)\\
  & = & \sum_{k\in\mathbb{Z}} \Big( \sum_{e(\alpha)\in\Omega_k \setminus \Omega_{k + 1}}| \varphi (\alpha) |^p \pi (\alpha) +\sum_{e(\alpha)\in\partial\Omega_k \cap \partial \Omega_{k+1}}\frac{| \varphi (\alpha) |^p \pi (\alpha)}{2^{p(m(\alpha)-k)}}\Big)\\
  & = & \sum_{\alpha}  \pi(\alpha) |\varphi(\alpha)|^p + \sum_\alpha \pi(\alpha) |\varphi(\alpha)|^p \sum_{k< m(\alpha)} \frac{ \chi_{\partial\Omega_{k}\cap \partial \Omega_{k+1}}(e(\alpha)) }{2^{p(m(\alpha)-k)}}\\
  & \leq & \frac{2^p}{2^p-1} \| \varphi \|_{\ell^p (E,\pi)}^p .
\end{eqnarray*}

\end{proof}


\begin{proof}[Proof of the equivalence of (i) and (iii) in Theorem \ref{main theorem}] Suppose that as always $\varphi $ is a positive function defined on the edges of the tree, then using the distribution function we write \begin{align*}
    \int_{\overline{T}} \III \varphi ^ p d\mu & = \int_{0}^\infty \mu(\{ x \in \overline{T}: \III\varphi (x) >  \lambda   \}) d\lambda^p \\
    & \leq (2^p-1) \sum_{k=-\infty}^{+\infty}2^{kp} \mu ( x \in \overline{T}: \III \varphi (x) > 2^k ) \\ 
    & \leq (2^p-1)[[\mu]]_c^p \sum_{k=-\infty}^{+\infty} 2^{kp} \Capac_{p,\pi}(\{ x\in \partial T : \III \varphi(x) > 2^k \}) \\ 
    & \leq 2^{p} [[\mu]]_c^p \norm{ \varphi }^p_{\ell^p(\pi)}.
\end{align*}

This concludes the proof of sufficiency.

To prove necessity let $\alpha_1, \dots, \alpha_n \in E$ and denote by $\varphi$ the equilibrium function associated to the set $\cup_{i=1}^n S(\alpha_i)$. Then, $I\varphi\geq 1$, $\Capac_{p,\pi}-$quasi everywhere and in particular $\mu-$a.e. on $\cup_{i=1}^n S(\alpha_i)$. It follows
\begin{equation*}
\begin{split}
     [\mu]^p\Capac_{p,\pi}(\cup_{i=1}^n S(\alpha_i))&=[\mu]^p\sum_{\alpha}\varphi(\alpha)^p \pi(\alpha)\\
     &\geq\int_{\overline{T}} (\III \varphi)^p d\mu\geq\sum_{i=1}^n\int_{S(\alpha_i)} (\III \varphi)^p d\mu\geq\sum_{i=1}^n\mu(S(\alpha_i)). 
\end{split} 
\end{equation*}
\end{proof}

\begin{problem}
Find the best constant in the inequality $[\mu] \leq 2^p [[\mu]]_c$ of Theorem \ref{main theorem}.
\end{problem}

\section{Mass--energy characterization: three different proofs}\label{sec:me}
In this section we will prove the mass--energy characterization for trace measures, that is, the equivalence of (i) and (ii) in Theorem \ref{main theorem}. We will give three different proofs, based on different techniques. The easiest proof only works for $p=2$, while the other two work for any $1<p<\infty$.

We remind that,
\begin{equation*}
 [[\mu]]^{\pp-1}=\sup_{\alpha\in E }\frac{\sum_{\beta\subseteq\alpha}\pi(\beta)^{1-\pp}\mu(S(\beta))^\pp}{\mu(S(\alpha))},
\end{equation*}
and we call it the \it energy-mass ratio\rm.

\subsection{Maximal function} This simple proof can be found, for instance, in \cite{ARSW2019}. It relies on the $L^p$ inequality for a suitable maximal function. If $\mu, \sigma \geq 0 $ are measures on $\overline{T}$, and $f$ a function on $\overline{T}$
\begin{equation}\label{DyadMax}
 M_\mu (fd\sigma) (x):=\max_{\alpha\in [o^*,x]} \frac{1}{\mu(S(\alpha))}\int_{S(\alpha)} |f| d\sigma.
\end{equation}
We simplify the notation by setting $M_\mu f=M_\mu (fd\mu)$. We have the following weak-$(1,1)$ estimate. 
\begin{theorem}\label{MaxThm}
 Suppose $\mu,\sigma\ge0$ are measures on $\overline{T}$ and $\psi$ a positive function on $\overline{T}$. Then,
 \begin{itemize}
  \item[(a)] $\displaystyle\lambda \sigma(x:\ M_\mu\psi(x)>\lambda)\le\int_{\overline{T}}\psi  M_\mu(d\sigma)d\mu$;
  \item[(b)] for $1<p<\infty$, $\displaystyle\int_{\overline{T}}(M_\mu\psi)^p d\sigma\le  (p^*)^p\int_{\overline{T}}\psi^p M_\mu(d\sigma)d\mu$.
 \end{itemize}
\end{theorem}
\begin{proof}
First we prove (a). Fix $\lambda >0$ and set
 $E(\lambda )=\{x\in \overline{T}:\ M_\mu\psi(x)> \lambda \}=\sqcup_{j}S(\alpha_j)$, where $\dashint_{S(\alpha_j)} \psi d\mu > \lambda$ and 
 $\dashint_{S(\beta)} \psi d\mu \le \lambda $ if $\beta \supset\alpha_j$. Then,
 \begin{eqnarray*}
  \sigma(E(\lambda ))&=&\sum_j\sigma(S(\alpha_j))=\sum_j\frac{\sigma(S(\alpha_j))}{\mu(S(\alpha_j))}\mu(S(\alpha_j))\crcr
  &\le&\frac{1}{\lambda}\sum_j\frac{\sigma(S(\alpha_j))}{\mu(S(\alpha_j))}\int_{S(\alpha_j)}\psi d\mu\le\frac{1}{\lambda}\sum_jM_\mu(d\sigma)(x_j)\int_{S(\alpha_j)}\psi d\mu \crcr
  &\le&\frac{1}{\lambda}\sum_j\int_{S(\alpha_j)}\psi M_\mu(d\sigma)d\mu\le\frac{1}{\lambda}\int_{\overline{T}}\psi M_\mu(d\sigma)d\mu,
 \end{eqnarray*}
which proves the weak estimate. Then a simple application of a variation of the classical Marcinkiewicz interpolation theorem \cite[Exercise 1.3.3]{Grafakos2008} gives us $(b)$ from $(a).$ 
\end{proof}

\begin{proof}[Proof of the equivalence of (i) and (ii) in Theorem \ref{main theorem}]
Assume the mass--energy condition holds, and define the measure $\sigma$ on $\overline{T}$ by setting  $\sigma(S(\alpha))=\sum_{\beta\subseteq\alpha}\mu(S(\beta))^\pp\pi(\beta)^{1-\pp}$.
The mass--energy condition  \eqref{EqMassEnergy} can be written as $\sigma(S(\alpha))\le[[\mu]]^{\pp}\mu(S(\alpha))$. Thus, for any positive function $\psi$ on $\overline{T}$,
\begin{eqnarray*}
\sum_{\alpha\in E}(I^*_\mu \psi)^\pp(\alpha)\pi(\alpha)^{1-\pp}&=& \sum_{\alpha\in E}\left(\dashint_{S(\alpha)}\psi d\mu\right)^\pp\sigma(e(\alpha))\crcr
&\le&\sum_{\alpha \in E} (M_\mu\psi)^\pp e((\alpha))\sigma(e(\alpha))\crcr
&\le & \int_{\overline{T}}(M_\mu\psi)^\pp (x)\sigma(x)\crcr
&\le & p^{p^*} \int_{\overline{T}}\psi^\pp M_\mu(d \sigma)d\mu\crcr
&\le& p^{p^*} [[\mu]]^{\pp-1}\int_{\overline{T}}\psi^\pp d\mu.
\end{eqnarray*}
So the dual Hardy's inequality \eqref{HardyInequalityDual}, equivalent to \eqref{HardyInequality} is obtained.

Conversely, suppose the dual Hardy's inequality \eqref{HardyInequalityDual} holds. By testing it on functions $\psi=\chi_{S(\alpha)}$, we obtain
\begin{equation*}
    \mu(S(\alpha))\sum_{\beta\supset\alpha}\pi(\beta)^{1-\pp} +\sum_{\beta\subseteq\alpha}\pi(\beta)^{1-\pp}\mu(S(\beta))^\pp\le [\mu]^{\pp-1}\mu(S(\alpha)),
\end{equation*}
which implies the mass--energy condition \eqref{EqMassEnergy}.
\end{proof}

\subsection{Monotone proof for p=2} This easiest proof only works for $p=2$, because it uses the $C^*$ identity for operators on Hilbert spaces. 

\begin{lemma}\label{TTstar}
 Let $T:H_1\to H_2$ be a bounded and linear operator between Hilbert spaces. Then, $\|T\|=\|T^*\|=\|TT^*\|^{1/2}$.
\end{lemma}

\begin{lemma}\label{LMonotone}
 Let $\mu,\nu\ge0$ be measures on $\overline{T}$, and suppose that $\mu(S(\alpha))\le\nu(S(\alpha))$ holds for all $\alpha$ in $E$; 
 and suppose that $f:\overline{T}\to\RR_+$ is monotone, $f(x)\ge f(y)$ if $x\subseteq y$.
 Then, $\int fd\mu\le\int fd\nu.$ 
\end{lemma}
\begin{proof} It suffices to prove the inequality of the corresponding distribution functions. For $t>0$,
 \begin{eqnarray*}
  \mu(x:\ f(x)>t)&=&\sum_j\mu(S(\alpha_j)) \text{ since $f$ is monotone,}\crcr
  &\le&\sum_j\nu(S(\alpha_j))=\nu(x:\ f(x)>t).
 \end{eqnarray*}
\end{proof}

\begin{proof}[Proof of the equivalence of (i) and (ii) in Theorem \ref{main theorem} for $p=2$]

Suppose initially that our tree is finite, but arbitrarily large. This assures that all relevant operators are bounded. If we manage to estimate the norm of the Hardy operator independently of the length of the tree we can pass to the infinite case by a simple limiting argument. 

For a $g:  E \to \RR_+ $, we first compute
 \begin{eqnarray*}
  \III^*\III g(\alpha)&=& \pi(\alpha)^{-1}\int_{S(\alpha)}\sum_{\beta\in[o^*,x]}g(\beta)d\mu(x)=\pi(\alpha)^{-1} \sum_{\beta}g(\beta)\mu(S(\beta)\cap S(\alpha))\crcr
  &=& \pi(\alpha)^{-1} \sum_{\beta \subseteq \alpha}g(\beta)\mu(S(\beta))+\pi(\alpha)^{-1}\mu(S(\alpha))\sum_{\beta\supseteq \alpha}g(\beta) \crcr
&=&\pi(\alpha)^{-1}\sum_{\beta\subseteq \alpha}g(\beta)\mu(S(\beta))+\pi(\alpha)^{-1}\mu(S(\alpha))\III g(e(\alpha))\crcr
&=& T_1g(\alpha)+T_2g(\alpha).
 \end{eqnarray*}
 
 Therefore, 
\begin{equation*}
	\norm{\III^*\III}_{\ell^2(\pi)} \leq \norm{T_1}_{\ell^2(\pi)}+\norm{T_2}_{\ell^2(\pi)}.
\end{equation*} 
The second norm can be computed in terms of the norm  $ \III $. Consider a measure $\rho $ on $T$ such that $\rho(e(\alpha)) = \pi^{-1}(\alpha)\mu(S(\alpha))$. The mass--energy condition allows us to apply Lemma \ref{LMonotone} to the measures $\rho$ and $\mu$. Therefore obtaining
\begin{eqnarray*}
\norm{T_2g}_{\ell^2(\pi)}^2&=& \sum_{\alpha}\pi(\alpha)^{-1} \mu(S(\alpha))^2\left(\III g(e(\alpha))\right)^2 = \int_T (\III g )^2 d\rho  \crcr
& \leq & [[\mu]] \int_T (\III g)^2 d\mu \crcr
& = & [[\mu]] \sum_{\alpha}\mu(\alpha)\left(\III g(e(\alpha))\right)^2 \leq [[\mu]] \norm{g}^2_{\ell^2(\pi)}\norm{\III}_{\ell^2(\pi)\to L^2(\mu)}^2.
\end{eqnarray*}

A standard calculation shows that $ T_2^*=T_1 $ (with respect to the inner product in $\ell^2(\pi)$), 
hence we have for free the estimate on the norm of $ T_1$. Putting everything together we get 
\begin{equation*}
	\norm{\III}^2_{\ell^2(\pi)\to L^2(\mu)} = \norm{\III\III^*}_{\ell^2(\pi)}\leq 2 [[\mu]]^{\frac{1}{2}} \norm{\III}_{\ell^2(\pi)\to L^2(\mu)}.
\end{equation*} Since $ \III $ is bounded because our tree is finite we can divide both sides of the inequality with its norm to get
\begin{equation*}
	\norm{\III}_{\ell^2(\pi)\to L^2(\mu)}^2 \leq 4 [[\mu]].
\end{equation*}
Notice that this proof gives the best constant. \end{proof}

\subsection{Bellman function} \label{sec:bellman}
In this section we provide a different proof, based on a Bellman function approach, of the fact that \eqref{EqMassEnergy} implies \eqref{HardyInequality}.
Let $T$ be a general rooted tree and denote by $|\cdot|$, as usual, the \textit{canonical} edge weight defined in Section \ref{subsec:Some potential theory}. Set $k_\beta=|\beta|/|\alpha|$ when $\beta\in s(\alpha)$.

The following is a tree version of the Weighted Dyadic Carleson Imbedding Theorem by Nazarov, Treil and Volberg \cite{NTV}. With respect to the standard dyadic case, this tree analogue presents some extra difficulty, due to the fact that the objects into play are here not martingales but only supermartingales.

\begin{theorem}[Carleson imbedding theorem for trees]\label{thm: CET}
Let $\sigma$ be a nonnegative weight on $E$ and $\lambda$ a measure on $\overline{T}$ satisfying

\begin{equation}\label{eq 1}
    \sum_{\beta\subseteq\alpha}\sigma(\beta)\frac{I_\lambda^*(\beta)^p}{|\beta|^p}\leq I_\lambda^*(\alpha), \quad \text{for every } \alpha\in E.
\end{equation}
Then
\begin{equation}\label{eq 2}
    \sum_{\alpha\in E}\sigma(\alpha)\frac{I_\lambda^*\varphi(\alpha)^p}{|\alpha|^p}\leq (p^*)^p \Vert \varphi\Vert_{\ell^p(\lambda)}^p, \quad \text{for every } \varphi:E\to \mathbb{R}_+.
\end{equation}


\end{theorem}

Letting $\mu$ be the measure defined by $\mu(S(\alpha))=I^*_\lambda(\alpha)/|\alpha|$, and switching to the $\pi-$dictionary under the usual correspondence $\pi(\alpha)=\sigma(\alpha)^{1-p}$, we see that \eqref{eq 1} is equivalent to the mass--energy condition \eqref{EqMassEnergy} and \eqref{eq 2} coincides with the dual Hardy's inequality \eqref{HardyInequalityDual}.

On the dyadic tree, a proof of the above theorem relying on a Bellman function method was first given in \cite{AHMV} for $p=2$, and later extended to every $1<p<\infty$ in \cite{Cavina2020}. In his paper it is proven that the Bellman function employed is \emph{the} Bellman function of a problem in stochastic optimal control. We give here a slightly adapted proof which works on every tree $T$. We also remark that the result remains true substituting the canonical weight $|\cdot|$ with a general weight $w$ fulfilling the so called \textit{flow condition}, that is, $\sum_{\beta\in s(\alpha)}w(\beta)=w(\alpha)$, for any $\alpha\in E$. Indeed, in the following proof the flow condition is the only property of the canonical weight which is used. We refer the reader to \cite{Chalmoukis2019}, \cite{lstv1} and \cite{lstv2} for some recent result concerning trees endowed with flow measures and flow weights.

\begin{proof}[Proof of Theorem \ref{thm: CET}]
We begin by observing that it is enough to show the result for nonnegative functions. For any edge $\alpha$ and any quadruple of nonegative real numbers $F,f,A,v$, define $\Omega_\alpha(F,f,A,v)$ to be the set of weights $\sigma$, measures $\lambda$ and functions $\varphi$ such that
\begin{equation*}
    \frac{1}{|\alpha|}I_\lambda^*(\varphi^p)(\alpha)=F, \quad \frac{1}{|\alpha|}I_\lambda^*\varphi(\alpha)=f, \quad \frac{1}{|\alpha|}\sum_{\beta\subseteq\alpha}\sigma(\beta)\frac{I_\lambda^*(\beta)^p}{|\beta|^p}=A, \quad \frac{1}{|\alpha|}I_\lambda^*(\alpha)=v.
\end{equation*}
In order for $\Omega_\alpha(F,f,A,v)$ not to be empty, it must be $f^p\leq F v^{p-1}$, the condition coming from Hölder inequality. Moreover, \eqref{eq 1} implies $A\leq v$. We denote by $\mathcal{D}$ the domain $\lbrace f^p\leq F v^{p-1}, \ A\leq v\rbrace\subseteq \mathbb{R}^4$, which is clearly convex, being the intersection of the half plane $\lbrace A\leq v\rbrace$ with the cylindroid having as a basis the convex set $\lbrace f^p\leq F v^{p-1}\rbrace\subseteq\mathbb{R}^3$. Define the \textit{Belmann function} $\mathbb{B}: \mathbb{R}^4\to\mathbb{R}_+$, 
\begin{equation*}
    \mathbb{B}(F,f,A,v)=\sup_{\Omega_\omega(F,f,A,v)}\frac{1}{|\alpha|}\sum_{\beta\subseteq \omega}\sigma(\beta)\frac{I_\lambda^*\varphi(\beta)^p}{|\beta|^p}.
\end{equation*}
We aim to prove that
\begin{equation*}
    \mathbb{B}(F,f,A,v)\leq |\omega|(p^*)^p F, \quad \text{for all } (F,f,A,v)\in \mathcal{D}.
\end{equation*}
Let $x=(F,f,A,v)$ be a point in $\mathcal{D}$ and fix arbitrarily $(\sigma,\lambda, \varphi)\in\Omega_\omega(x)$. For each $\alpha\in E$, let $x_\alpha=(F_\alpha,f_\alpha,A_\alpha,v_\alpha)\in\mathbb{R}^4$ be the unique point such that $(\sigma,\lambda, \varphi)\in\Omega_\alpha(x_\alpha)$. In particular $x_\omega=x$, and it is clear that $x_\alpha\in \mathcal{D}$ for each $\alpha$. Moreover, the additivity of $I^*$ gives the following relations
\begin{equation*}
\begin{split}
    &F_\alpha=\frac{1}{|\alpha|}\varphi^p(\alpha)\lambda(\alpha)+\sum_{\beta\in s(\alpha)}k_\beta F_\beta, \quad f_\alpha=\frac{1}{|\alpha|}\varphi(\alpha)\lambda(\alpha)+\sum_{\beta\in s(\alpha)}k_\beta f_\beta,\\ &A_\alpha=\sigma(\alpha)\frac{I_\lambda^*(\alpha)^p}{|\alpha|^{p+1}}+\sum_{\beta\in s(\alpha)}k_\beta A_\beta, \quad v_\alpha=\frac{\lambda(\alpha)}{|\alpha|}+\sum_{\beta\in s(\alpha)}k_\beta v_\beta.
\end{split}
\end{equation*}
By denoting $b_\alpha^p=|\alpha|^{-1}\varphi^p(\alpha)$, $c_\alpha=\sigma(\alpha)|\alpha|^{-(p+1)}I_\lambda^*(\alpha)^p$ and $a_\alpha^{p^*}=|\alpha|^{-1}\lambda(\alpha)$, and defining the point $y_\alpha=(b_\alpha^p, a_\alpha b_\alpha, c_\alpha, a_\alpha^{p^*})$, the above relations can be rewritten as
\begin{equation}\label{eq: recursive family-short}
      x_\alpha=y_\alpha+\sum_{\beta\in s(\alpha)}k_\beta x_\beta, \quad \alpha\in E.
\end{equation}
Now, suppose we can design a concrete function, $\mathcal{B}:\mathbb{R}^4\to\mathbb{R}$ such that $(i) \ \mathcal{B}(F,f,A,v)\leq (p^*)^pF$ on $\mathcal{D}$, and satisfying
\begin{equation*}
    (ii) \  |\alpha|\mathcal{B}(x_\alpha)-\sum_{\beta\in s(\alpha)}|\beta|\mathcal{B}(x_\beta)\geq \sigma(\alpha)f_\alpha^p, \quad \text{for every } \alpha\in E.
\end{equation*}
Then, summing over $\alpha\in E$ both sides of the inequality and exploiting the telescopic structure of the summand we obtain
\begin{equation}
     \sum_{\alpha\in E} \sigma(\alpha)\frac{1}{|\alpha|}I^*_\lambda(\varphi^p)(\alpha)=\sum_{\alpha\in E} \sigma(\alpha)f_\alpha^p\leq |\omega|\mathcal{B}(x)\leq |\omega|(p^*)^p F,
\end{equation}
from which follows the thesis,
\begin{equation*}
    \mathbb{B}(F,f,A,v)\leq |\omega|(p^*)^p F, \quad \text{for all } (F,f,A,v)\in \mathcal{D}.
\end{equation*}
We now claim that the function
\begin{equation*}
    \mathcal{B}(F,f,A,v)=(p^*)^p\Big(F-\Big(\frac{p-1}{A+(p-1)v}\Big)^{p-1}f^p\Big),
\end{equation*}
fulfils the desired properties. The construction of a Bellman's functions is a delicate matter. The interested reader can find more information and examples in \cite{Burkholder1984, NTV2000, banuelos2010donald}.   It is immediate that $(i)$ holds on $\mathcal{D}$. For any chosen $x=(F,f,A,v)\in \mathcal{D}$, let $x_\alpha$ be the associated family of points solving \eqref{eq: recursive family-short}. Then also the points $x_\alpha^*=x_\alpha-(0,0,c_\alpha,0)$ and $x_\alpha^{**}=x_\alpha-y_\alpha$ belong to the convex domain $\mathcal{D}$. Since $\mathcal{B}$ is clearly concave in the third variable, we have
\begin{equation*}
    \mathcal{B}(x_\alpha)-\mathcal{B}(x_\alpha^*)\geq c_\alpha\frac{\partial \mathcal{B}}{\partial A}(x_\alpha)=c_\alpha\Big(\frac{pf_\alpha}{A_\alpha+(p-1)v_\alpha}\Big)^p\geq c_\alpha\Big(\frac{f_\alpha}{v_\alpha}\Big)^p,
\end{equation*}
the last inequality following from the domain constraint $A_\alpha\leq v_\alpha$. Indeed, $\mathcal{B}$ is also concave as a function of four variables on the convex set $\mathcal{D}$, as one can verify by checking that the Hessian matrix of $(F,f,A,v)\mapsto \mathcal{B}(F,f,A,v)$ is positive semidefinite on $\mathcal{D}$. Hence, we have
\begin{equation*}
     \mathcal{B}(x_\alpha^*)-\mathcal{B}(x_\alpha^{**})\geq \frac{\partial \mathcal{B}}{\partial F}(x_\alpha^*)b_\alpha+\frac{\partial \mathcal{B}}{\partial f}(x_\alpha^*)a_\alpha b_\alpha+\frac{\partial \mathcal{B}}{\partial v}(x_\alpha^*)a_\alpha^{p^*}\geq 0, 
     \end{equation*}
where the last inequality can be derived by direct calculations. Putting the pieces together we obtain
\begin{equation}\label{eq: main special single step-short}
    c_\alpha\Big(\frac{f_\alpha}{v_\alpha}\Big)^p\leq \mathcal{B}(x_\alpha)-\mathcal{B}(x_\alpha^*)=\mathcal{B}(x_\alpha)-\mathcal{B}(x_\alpha^{**})+\mathcal{B}(x_\alpha^{**})-\mathcal{B}(x_\alpha^*)\leq \mathcal{B}(x_\alpha)-\mathcal{B}(x_\alpha^{**}).
\end{equation}
Exploiting the concavity of $\mathcal{B}$,
\begin{equation*}
    \mathcal{B}(x_\alpha^{**})=\mathcal{B}\Big(\sum_{\beta\in s(\alpha)}k_\beta x_\beta\Big)\geq \sum_{\beta\in s(\alpha)}k_\beta\mathcal{B}(x_\beta),
\end{equation*}
which, by means of \eqref{eq: main special single step-short} yields to
\begin{equation}
    c_\alpha\Big(\frac{f_\alpha}{v_\alpha}\Big)^p\leq \mathcal{B}(x_\alpha)-\sum_{\beta\in s(\alpha)}k_\beta\mathcal{B}(x_\beta).
\end{equation}
It is easily seen that $c_\alpha=v_\alpha^p\sigma(\alpha)|\alpha|^{-1}$, which substituted above gives $(ii)$.
\end{proof}

It is clear that, a posteriori, the mass--energy and the isocapacitary conditions are equivalent, being both equivalent to \eqref{HardyInequality}. However, it is tempting to look for an a priori argument for the equivalence of these geometric conditions which does not pass through the boundedness of the Hardy operator. A direct proof that \eqref{EqMassEnergy} implies \eqref{eq:isocapacitary}, for a family of weights including $\pi = 1$, is in \cite{ARS2008} ,
where it is also directly proven, for $p = 2$ and $\pi = 1$, that $\nu\leq \mu$ implies that $[[\nu]] \leq 2 [[\mu]]$.

\begin{problem}
Find a proof of the equivalence between \eqref{EqMassEnergy} and \eqref{eq:isocapacitary}, which works for every couple $\pi,\mu$ and does not require the boundedness of the Hardy operator.
\end{problem}

\section{A reverse H\"older inequality}\label{sec:holder}

In the particular case that $\pi \equiv 1 $ and $p=2$ the mass--energy condition can be rewritten in an interesting way as a consequence of the following calculation. 
\begin{align*} \sum_{\beta \subseteq \alpha } \mu(S(\beta))^2  & = \sum_{\beta \subseteq \alpha } \int_{S(\beta)} 
\int_{S(\beta)} d\mu(x) d\mu(y) \\
& =   \int_{S(\alpha)} \int_{S(\alpha)} \Big( \sum_{\beta \supseteq x \wedge y} 1 \Big) d\mu(x)d\mu(y) \\
& = \int_{S(\alpha)} \int_{S(\alpha)} d(x\wedge y) d\mu(x)d\mu(y).
\end{align*}
Therefore, the mass--energy condition can be expressed as 
 \[ \sup_{\alpha \in E } \dashint_{S(\alpha)} \int_{S(\alpha)} d(x\wedge y) d\mu(x) d\mu(y) < + \infty.  \]
 Notice also that 
 \[ \int_{S(\alpha)} d(x\wedge y) d\mu(x) = \III\III_\mu^* (\chi_{S(\alpha)}). \]
The following variation on the above condition,
\[  \sup_{\alpha \in E} \dashint_{S(\alpha)} \Big( \int_{S(\alpha)} d(x\wedge y) d\mu(x) \Big)^s d\mu(y): = [[\mu]]_s < \infty, \tag{s-Testing} \]
 is clearly stronger than the mass--energy condition for $s>1$ due to H\"older's inequality. 
 The surprising result is that in fact the conditions are equivalent, and the corresponding quantities are comparable. This result is in the spirit of the John-Nirenberg reverse H\"older inequality for BMO functions.
 
 \begin{theorem}\label{CalderonZygmound} For all measures $\mu$, and $s > 1,$
  \[ \sup_{\alpha \in E} \dashint_{S(\alpha)} \Big( \int_{S(\alpha)} d(x\wedge y) d\mu(x) \Big)^s d\mu(y) \leq C_s [[ \mu ]]^s.\]
 \end{theorem}
 
 In an implicit form this result is contained in the work of Tchoundja \cite{Tchoundja2008}. The proof of the above theorem is based on a Calderón–Zygmund type theorem for the operator $ \III\III_\mu^* : = T_\mu $. More precisely, 
\begin{theorem}\label{reverse holder} Suppose that the operator 
\[ T_\mu : L^2(\mu) \to L^2(\mu) \] is bounded. Then for any $s \in (1,+\infty)$ the operator 
\[ T_\mu : L^s(\mu ) \to L^s(\mu) \]
is bounded.  Furthermore,  $ \norm{T_\mu}_{L^s(\mu)} \leq C_s \norm{T_\mu}_{L^2(\mu)}. $
\end{theorem}

Since the underlying measure $\mu $ is not necessarily doubling this theorem can be seen as a special case of \cite[Theorem 1.1]{Nazarov1998}. Here we shall give a direct proof which also provides better quantitative estimates of the constants involved based on a good-$\lambda$ inequality as in \cite{Verdera2000, Tchoundja2008}.

\begin{lemma}[Good-$\lambda$ inequality]\label{lem:good-lambda}
 Let $\mu$ a trace measure on $\overline{T}$. Then for every $\eta>0$, there exists $\gamma(\eta)>0$ such that for any non negative function $f$ on $\overline{T}$,
 \[ \mu\{x\in \overline{T}: T_\mu f(x) > (1+\eta) \lambda \,\,\, \text{and} \,\,\, M_\mu(f^2)(x) \leq \gamma^2\lambda ^2 \}  \leq \frac{1}{2}\mu \{ x\in \overline{T}: T_\mu f(x) > \lambda  \}. \]
\end{lemma}

\begin{proof}
Notice that the set $\{ T_\mu f>\lambda  \}$ is a {\it stopping time}. In other words it can be written as a disjoint union of tent regions,
\begin{equation*}
    \{T_\mu f > \lambda  \} = \bigcup_{i=1}^\infty S(\alpha_i).
\end{equation*}
It is therefore sufficient to prove that for all $\alpha_i$ we have 
 \[ \mu\{x\in S(\alpha_i): T_\mu f(x) > (1+\eta) \lambda \,\,\, \text{and} \,\,\, M_\mu(f^2)(x) \leq \gamma^2\lambda ^2 \}  \leq \frac{1}{2}\mu(S(\alpha_i)). \]
 So for the rest of the proof we work on a fixed $S(\alpha_i)$ which we denote by $S(\alpha)$ to avoid an overload of notation. Let $f_1=f\chi_{S(\alpha)}$ and $f_2=f-f_1$. For $x\in S(\alpha)$
 \begin{align*} T_\mu f_2 (x) & = \int_{\overline{T}\setminus S(\alpha)}d(x\wedge y) f(y) d\mu(y) \\
  & = \int_{\overline{T}\setminus S(\alpha)} d(b(\alpha)\wedge y) f(y)d\mu(y) \leq T_\mu f(b(\alpha))  \leq \lambda .
  \end{align*}
 Because $ b(\alpha) \not\in S(\alpha). $
  Therefore $ T_\mu f (x) \leq T_\mu f_1 (x) + \lambda $, which implies that
  \begin{align*}  \mu\{x\in S(\alpha): T_\mu f(x) > (1+\eta) \lambda \,\,\,&  \text{and} \,\,\, M_\mu(f^2)(x) \leq \gamma^2\lambda ^2 \} \\
   & \leq  \mu\{x\in S(\alpha): T_\mu f_1(x) > \eta  \lambda  \}  \\
   & \leq \frac{1}{\eta^2 \lambda ^2} \int_{S(\alpha)} (T_\mu f_1) ^2 d\mu \\ 
   & \leq \frac{\norm{T_\mu}^2_{L^2(\mu)}}{\eta^2\lambda^2 }\int_{S(\alpha)}f^2 d\mu \\ 
   & \leq  \frac{\norm{T_\mu}^2_{L^2(\mu)} \mu(S(\alpha))}{\eta^2\lambda^2 } M_\mu(f^2)(e(\alpha)) \\
   & \leq  \frac{\norm{T_\mu}^2_{L^2(\mu)}\gamma ^2 }{\eta^2  }
 \mu(S(\alpha)).  \end{align*}
  Where we assume without loss of generality that  $M_\mu(f^2)(e(\alpha)) \leq \gamma^2 \lambda ^2$, otherwise the left hand side is zero. It suffices therefore to choose 
  \[  \gamma = \frac{\eta }{\sqrt{2} \norm{T_\mu}_{L^2(\mu)} }.\]
\end{proof}

\begin{proof}[Proof of Theorem \ref{CalderonZygmound}]
Since the operator $T_\mu$ is self adjoint it suffices to prove that $L^2(\mu)$ boundedness implies $ L^s(\mu)$ boundedness for all $s>2$. Let $s>2 $ and $f\in L^s(\overline{T},\mu)$. Exploiting Lemma \ref{lem:good-lambda} and Theorem \ref{MaxThm}, we get
\begin{align*}
    \int_{\overline{T}} (T_\mu f)^s d\mu & = \int_0^\infty \mu(\{ T_\mu f > \lambda \}) d\lambda ^s \\ 
    & = (1+\eta)^s  \int_0^\infty \mu(\{ T_\mu f > (1+\eta) \lambda \}) d\lambda ^s \\
    &  \leq (1+\eta)^s  \int_0^\infty \mu(\{ T_\mu f > (1+\eta) \lambda \,\,\, \text{and} \,\,\, M_\mu(f^2) \leq \gamma^2 \eta^2 \}) d\lambda ^s \\ & \,\,\,\,\,\,\,\,\,  + (1+\eta)^s  \int_0^\infty \mu(\{ M_\mu(f^2) > \gamma^2\lambda^2  \}) d\lambda^s \\ 
   &  \leq \frac{(1+\eta)^s}{2}  \int_0^\infty \mu(\{ T_\mu f > \lambda \}) d\lambda ^s + \frac{(1+\eta)^s}{\gamma^s} \int_{\overline{T}} M_\mu(f^2)^\frac{s}{2} d\mu \\
   & \leq \frac{(1+\eta)^s}{2} \int_{\overline{T}} (T_\mu f)^s d\mu + \frac{2^{s/2}s(1+\eta)^s}{(s-2)\gamma^s} \int_{\overline{T}} f^s d\mu,
\end{align*}
which proves the thesis if $\eta$ is chosen small. In particular 
\[ \norm{T_\mu}_{L^s(\mu)} \leq C_s \norm{T_\mu}_{L^2(\mu)}. \] \end{proof}

The reverse H\"older inequality is now a corollary of the above theorem.

\begin{proof}[Proof of Theorem \ref{reverse holder}]
Suppose that $\mu$ satisfies the mass--energy condition. Then, 
\[ \norm{T_\mu}_{L^2(\mu)} = \norm{\III \III^*_\mu}_{L^2(\mu)} = \norm{\III}^2_{\ell^2 \to L^2(\mu)} \leq 4 [[\mu]]. \]

On the other hand,  
\[ \norm{T_\mu}_{L^s(\mu)}^s \geq \frac{\norm{T_\mu(\chi_{S(\alpha))})}_{L^s(\mu)}^s}{\mu(S(\alpha))} \geq \dashint_{S(\alpha)} \III\III^*_\mu(\chi(S(a
))^s d\mu, \]
and the result follows from Theorem \ref{CalderonZygmound}.
\end{proof}

\section{The inequality of Muckenhoupt and Wheeden, and Wolff}\label{sec:MW}

In this section we only consider only the case when $T$ is a homogeneous tree. We recall that if each vertex of $T$ has $q+1$ neighbors, then $|\alpha|=q^{-d(\alpha)}$, for each edge $\alpha$. Let $0 < s < 1$ and $1 < p^* < \infty$. For any measure $\mu$ on $\overline{T}$ we trivially have
\[ \left( \sup_{\alpha \supset x} \frac{\mu (S (\alpha))}{|\alpha|^s}
   \right)^{p^*} \leq \sum_{\alpha \supset x} \left( \frac{\mu (S
   (\alpha))}{|\alpha|^s} \right)^{p^*} \leq \left( \sum_{\alpha
   \supset x} \frac{\mu (S (\alpha))}{|\alpha|^s} \right)^{p^*} \quad \forall x\in \partial T. \]
The inequality of Muckenhoupt and Wheeden \cite[Theorem 1]{Muckenhoupt1974}, \eqref{eq:MW} in the sequel, says
that the chain of inequalities can be reversed, on average.\footnote{In fact the the full Muckenhoupt--Wheeden inequality in the Euclidean setting, applies to more general situations when the underling measure is only $\mathcal{A}_\infty$ equivalent to the Lebesgue measure $dx$.}

\begin{theorem}
  For any measure $\mu$ on $\partial T$, $p^*\geq 1$, and $0<s<1$, there is a constant $C = C (p^*s)$ such that
  \begin{equation}\label{eq:MW}\tag{MW}
     \frac{1}{C} \int_{\partial T} \left( \sum_{\alpha \supset x} \frac{\mu (S
     (\alpha))}{|\alpha|^s} \right)^{p^*} d x \leq \int_{\partial
     T} \sum_{\alpha \supset x} \left( \frac{\mu (S (\alpha))}{|\alpha|^s}
     \right)^{p^*} d x \leq C \int_{\partial T} \left( \sup_{\alpha \supset x}
     \frac{\mu (S (\alpha))}{|\alpha|^s} \right)^{p^*} d x. 
  \end{equation}
\end{theorem}
As usual, $d x$ here is the Lebesgue measure for which $\int_{S (\alpha)} d x |\alpha|$. A first consequence of the \eqref{eq:MW} inequality is that we have a one parameter of
seemingly different conditions characterizing $\mu$'s for which the Hardy
inequality holds, provided that the weight $\pi$ has the special form $\pi(\alpha) = |\alpha|^{\frac{p^*s-1}{1-p^*}}$. Indeed, the central term in \eqref{eq:MW} can be written as an energy,
\begin{equation*}
    \begin{split}
        \int_{\partial T} \sum_{\alpha \supset x} \left( \frac{\mu (S (\alpha))}{|\alpha|^s} \right)^{p^*} d x &=  \sum_{\alpha}\frac{\mu (S (\alpha))^{p^*}}{|\alpha|^{- p^*s}} \int_{\partial S (\alpha)} d x \\
   &= \sum_{\alpha}
   \mu (S (\alpha))^{p^*} |\alpha|^{(1-p^* s)}=\mathcal{E}_{p,\pi}(\mu) ,
    \end{split}
\end{equation*}
and the \eqref{eq:MW} gives $\mathcal{E}_{p,\pi}(\mu)\approx \int_{\partial T} \left( \sum_{\alpha \supset x} \left( \frac{\mu (S
     (\alpha))}{|\alpha|^s}\right)^q \right)^{p^*} d x$, for all $q\geq 1$.

\begin{proposition}[Wolff's inequality on the tree]\label{prop:Wolff's inequality on the tree}
Let $\mu$ be a non-negative Borel measure on $\partial T$. Then for any $p^*\geq 1$ and $0<s<1$ one has
\begin{equation}\label{e:40}
\int_{\partial T }\left(\sum_{\alpha \supset x}\frac{\mu(S(\alpha) )}{|\alpha|^s}\right)^{p^*}\,dx \lesssim \int_{\partial T}\sum_{\alpha \supset x}\frac{(\mu(S(\alpha) ))^{p^*}}{|\alpha|^{s p^*}}\,dx.
\end{equation}
\end{proposition}

Since the particular choice of $q$ plays no role, from now on, to keep the notation lighter, we fix the homogeneity of the tree $T$ setting $q=2$., i.e., we put ourselves back in the realm of the classical dyadic Hardy's inequality \eqref{HardyDyadicIntervals}. Given an edge $\alpha\in E$ we denote by $\alpha^+$ and $\alpha^-$ its two children edges. In this setting, Proposition \ref{prop:Wolff's inequality on the tree} follows from the slightly more general statement below. The function $\varphi: T \rightarrow \mathbb{R}_+$ is a logarithmic supermartingale with the drift $d>0$, if for every edge $\alpha$ one has
\begin{equation}\label{e:100}
\frac12\left(\log\varphi(\alpha^+) + \log\varphi(\alpha^-)\right) \leq \log\varphi(\alpha) - d.
\end{equation}
\begin{proposition}\label{p:2}
Assume $\varphi$ is a logarithmic supermartingale with the drift $d>0$, and that its jumps are bounded from above,
\begin{equation}\label{eq:bounded jump}
 \max(\varphi(\alpha^+),\varphi(\alpha^-)) \leq C\varphi(\alpha), \quad \alpha \in E,   
\end{equation}
for some constant $C>0$. Then, for any $p^*\geq 1$ one has
\begin{equation}\label{e:101}
\int_{\partial T}\left(\sum_{\alpha \supset  x}\varphi(\alpha)\right)^{p^*}\,dx \leq C_1(p^*,d)\int_{\partial T}\sum_{\alpha \supset x}\varphi^{p^*}(\alpha)\,dx \leq C_2(p^*,d)\int_{\partial T}\sup_{\alpha \supset x}\varphi^{p^*}(\alpha)\,dx.
\end{equation}
\end{proposition}

\begin{proof}[Proof of Proposition \ref{prop:Wolff's inequality on the tree}]
One only needs to observe that $\varphi(\alpha):= \mu(S(\alpha))/|\alpha|^s$, $\alpha\in E$, defines a logarithmic supermartingale with the drift $\log2\cdot(1-s)$ and bounded jumps. Indeed, given any edge $\alpha$ in $T$, we clearly have $|\alpha^{\pm}|^s = 2^{-s}|\alpha|^s$, hence, since $\mu(S(\alpha)) = \mu(S(\alpha^+)) + \mu(S(\alpha^-))$, we see that
\begin{equation}\notag
\begin{split}
&\frac{\mu(S(\alpha^+))}{|\alpha^+|^s}\cdot \frac{\mu(S(\alpha^-))}{|\alpha^-|^s} = 2^{2s}|\alpha|^{-2s}\mu(S(\alpha^+))\mu(S(\alpha^-))\\ & \leq 
|\alpha|^{-2s}2^{2(s-1)}(\mu(S(\alpha^+))+ \mu(S(\alpha^-)))^2 = 2^{2(s-1)}\left(\frac{\mu(S(\alpha))}{|\alpha|^s}\right)^2.
\end{split}
\end{equation}
The logarithmic supermartingale property follows immediately. On the other hand,
\[
\frac{\mu(S(\alpha^{\pm}))}{|\alpha^{\pm}|^s} \leq 2^s\frac{\mu(S(\alpha))}{|\alpha|^s},
\]
so the jumps of $\varphi$ are clearly bounded from above.
\end{proof}

In order to prove Proposition \ref{p:2} we will need the following lemma, of which we postpone the proof.
\begin{lemma}[Wolff’s  lemma]\label{Wolff's lemma}
  Fix $\delta>0$ and $N >1/2$. Let $\varphi$ be a logarithmic supermartingale with positive drift $d>0$ satisfying \eqref{eq:bounded jump}. Then for any edge $\alpha_0$ in $T$ the following inequality holds
  \begin{equation}\label{e:51}
      \sum_{\alpha\subseteq \alpha_0}\varphi^{\delta+N}(\alpha)|\alpha|\gtrsim  \sum_{\alpha\subseteq \alpha_0}\varphi^\delta(\alpha)\sum_{\beta\subseteq \alpha}\varphi^N(\beta)|\beta|.
  \end{equation}
\end{lemma}
\begin{proof}[Proof of Proposition \ref{p:2}]
We only show the left inequality in \eqref{e:101}.
Let us introduce following notations: write $[p^*]$ and $\{p^*\}$ for the integer and the decimal part of $p$ and set
\begin{equation*}
    \begin{split}
       &p^* = [p^*] + \{p^*\},\\
    &r = \frac{1 + \{p^*\}}{2},\\
    &Q = [p^*] - 1. 
    \end{split}
\end{equation*}
Note that $p^* = Q + r + r$. Then we have
\begin{equation}\label{e:102}
\begin{split}
&\int_{\partial T}\left(\sum_{\alpha \supset x}\varphi(\alpha)\right)^q\,dx\\
&\leq\int_{\partial T}\left(\sum_{\alpha \supset x}\varphi(\alpha)\right)^{Q}\left(\sum_{\alpha\supset x}\varphi^{p}(\alpha)\right) \left(\sum_{\alpha\supset x}\varphi^{p}(\alpha)\right)\,dx  \\
&=\int_{\partial T}\left(\sum_{\alpha_1,\dots, \alpha_{Q}, \alpha_{Q+1}, \alpha_{Q+2} \supset x}\varphi(\alpha_1)\dots\varphi(\alpha_{Q})\varphi^{p}(\alpha_{Q+1})\varphi^{p}(\alpha_{Q+2})\right)\,dx\\
&=\sum_{\pi\in S_{Q+2}}\sum_{\alpha_1\supset\dots\supset \alpha_{Q+1}\supset \alpha_{Q+2}}\varphi^{p_{\pi(1)}}(\alpha_{1})\dots\varphi^{p_{\pi(Q+2)}}(\alpha_{Q+2})|\alpha_{Q+2}|,
\end{split}
\end{equation}
where $S_{Q+2}$ is the symmetric group of all permutations of $\{1,\dots,Q+2\}$ and $p_j = 1,\; 1\leq j\leq Q$, $p_{Q+1} = p_{Q+2} = p$.\\
The next step is to use Wolff's lemma: given a permutation $\pi\in S_{Q+2}$ we apply \eqref{e:51} repeatedly to \eqref{e:102}, obtaining 

\begin{align*}\notag
&\sum_{\alpha_1\supset\dots\supset \alpha_{Q+1}\supset \alpha_{Q+2}}\varphi^{p_{\pi(1)}}(\alpha_{1})\dots\varphi^{p_{\pi(Q+2)}}(\alpha_{Q+2})|\alpha_{Q+2}|\\
&=\sum_{\alpha_1\supset\dots\supset \alpha_{Q}}\varphi^{p_{\pi(1)}}(\alpha_{1})\dots\varphi^{p_{\pi(Q)}}(\alpha_{Q})\left(\sum_{\alpha_Q\supset \alpha_{Q+1} \supset \alpha_{Q+2}}\varphi^{p_{\pi(Q+1)}}(\alpha_{Q+1})\varphi^{p_{\pi(Q+2)}}(\alpha_{Q+2})|\alpha_{Q+2}|\right)\\
&\lesssim\sum_{\alpha_1\supset\dots\supset \alpha_{Q}}\varphi^{p_{\pi(1)}}(\alpha_{1})\dots\varphi^{p_{\pi(Q)}}(\alpha_{Q})\left(\sum_{\alpha_Q\supset \alpha_{Q+1} }\varphi^{p_{\pi(Q+1)}+p_{\pi(Q+2)}}(\alpha_{Q+1})|\alpha_{Q+1}|\right)\\
&=\sum_{\alpha_1\supset\dots\supset \alpha_{Q-1}}\varphi^{p_{\pi(1)}}(\alpha_{1})\dots\varphi^{p_{\pi(Q-1)}}(\alpha_{Q-1})\times\\
&\times\left(\sum_{\alpha_{Q-1}\supset \alpha_{Q}\supset \alpha_{Q+1} }\varphi^{p_{\pi(Q)}}(\alpha_{Q})\varphi^{p_{\pi(Q+1)}+p_{\pi(Q+2)}}(\alpha_{Q+1})|\alpha_{Q+1}|\right)\\
&\lesssim\dots \lesssim \sum_{\alpha_1}\varphi^{p_{\pi(1)}+\dots+p_{\pi_{Q+2}}}(\alpha_1)|\alpha_1| =  \int_{\partial T}\sum_{\alpha \supset x}\varphi^q(\alpha)\,dx.
\end{align*}
Summing over all $\pi\in S_{Q+2}$ we obtain first half of \eqref{e:101}.
\end{proof}

We now prove Wolff's Lemma. The proof is based on a careful analysis on slow and fast growing geodesics.

\begin{proof}[Proof of Lemma \ref{Wolff's lemma}]
Without any loss of generality we may assume that $N=1$ (since the proof works all the same for every $N$) and $\alpha_0 = \omega$, the root edge.


What we are going to do next is to fix an edge $\alpha$ and look at the possible growth rate of $\varphi(\beta)$ for $\beta\subseteq\alpha$. The idea is that, if $\varphi$ does not grow too fast in this region, then one could expect for the second sum on the right hand side of \eqref{e:51} to be dominated by the value at the starting point,
\[
\sum_{\beta\subseteq\alpha}\varphi(\beta)|\beta|\lesssim \varphi(\alpha)|\alpha|.
\]
On the other hand, if $\varphi$ grows very (exponentially) fast, then we write the right-hand side of \eqref{e:51} as $\sum_{\beta}\varphi(\beta)|\beta|\sum_{\alpha_0 \supseteq\alpha \supseteq \beta}\varphi^{\delta}(\alpha)$, and expect the second sum to be estimated by the value at $\beta$.

Given $\alpha\in E $ and $k\geq0$, let 
\begin{equation}\notag
A(\alpha,k) := \{\beta\subseteq\alpha: \ d(\beta,\alpha) = k, \  \frac{\varphi(\beta)}{\varphi(\alpha)} \leq (1+r)^{k}\},
\end{equation}
to be the set of slowly growing successors (here $r = r(d,\delta)<10^{-2}$ is some small constant to be chosen later), let also $A(\alpha) = \bigcup_{k\geq0}A(\alpha,k).$ We have
\begin{equation}\notag
\begin{split}
\sum_{\alpha\in E}\varphi^{\delta}(\alpha)\sum_{\beta\subseteq\alpha}\varphi(\beta)|\beta| = \sum_{\alpha\in E}\varphi^{\delta}(\alpha)\sum_{\beta\in A(\alpha)}\varphi(\beta)|\beta| +\sum_{\alpha\in E}\varphi^{\delta}(\alpha)\sum_{\beta\subseteq\alpha,\beta\notin A(\alpha)}\varphi(\beta)|\beta|.
\end{split}
\end{equation}
We start by estimating the second term,
\begin{equation}\label{e:62}
\begin{split}
\sum_{\alpha\in E}\varphi^{\delta}(\alpha)\sum_{\beta\subseteq\alpha,\beta\notin A(\alpha)}\varphi(\beta)|\beta| &= \sum_{\beta\in E}\varphi(\beta)|\beta|\sum_{\alpha \supseteq \beta, \beta\notin A(\alpha)}\varphi^{\delta}(\alpha) \\
&<\sum_{\beta\in E}\varphi(\beta)|\beta|\sum_{\alpha \supseteq \beta, \beta\notin A(\alpha)}\varphi^{\delta}(\beta)(1+r)^{-\delta d(\alpha,\beta)}\\
&\leq C(r)\sum_{\beta\in E}\varphi^{\delta+1}(\beta)|\beta|.
\end{split}
\end{equation}
To deal with the first term we let
\[
B(\alpha,k) := \{\beta\subseteq\alpha,\; d(\beta,\alpha) = k: \frac{\varphi(\beta)}{\varphi(\alpha)} \geq (1-r)^k\},\quad \alpha\in E,\; k\geq0,
\]
and, as before, $B(\alpha) = \bigcup_{k\geq0}B(\alpha,k)$. The function $\varphi$ decays exponentially outside of $B(\alpha)$, in particular
\begin{equation}\label{e:63}
\begin{split}
&\sum_{\alpha\in E}\varphi^{\delta}(\alpha)\sum_{\beta\subseteq \alpha, \ \beta\notin B(\alpha)}\varphi(\beta)|\beta| < \sum_{\alpha\in E}\sum_{\beta\subseteq \alpha, \ \beta\notin B(\alpha)}(1-r)^{d(\alpha,\beta)}\varphi^{1+\delta}(\alpha)|\beta|\\
&\leq C(r)\sum_{\alpha\in E}\varphi^{1+\delta}(\alpha)|\alpha|.
\end{split}
\end{equation}
So far, we took care of two types of behaviour of $\varphi$: points of very fast growth (i.e., $\beta\notin A(\alpha)$), and points of very fast decay ($\beta\notin B(\alpha)$). Now we consider the points $\beta\subseteq\alpha$ where $\varphi(\beta)$ is roughly comparable to $\varphi(\alpha)$. It turns out that these points are very rare in the successor set of $\alpha$. More precisely we show that for every $\alpha\in E$ and $k\geq 0$
\begin{equation}\label{e:71}
|B(\alpha,k)| = \sharp\{\beta\in B(\alpha,k)\} \leq C(r)2^{\frac{k}{2}}.
\end{equation}The reason for this is that by the multiplicative property \eqref{e:100} the function $\varphi$ decays exponentially on (geometric) average, and its pointwise growth rate is bounded from above.

Let $\Phi = \log\varphi$. By \eqref{eq:bounded jump},
\[
\max\left(\frac{\varphi(\beta^+)}{\varphi(\beta)}, \frac{\varphi(\beta^-)}{\varphi(\beta)}\right) \lesssim 1, 
\]
hence
\[
\max\left(\Phi(\beta^+)-\Phi(\beta), \Phi(\beta^-) - \Phi(\beta)\right) \lesssim 1.
\]
On the other hand (remind that $r< 10^{-2}$), if $\beta\in B(\alpha,k)$, then we get 
\[
\Phi(\beta) - \Phi(\alpha) \geq -kr = -d(\alpha,\beta)r.
\]
Choose $r< \min (1/100, d/10)$.
The inequality \eqref{e:71} now follows from the following lemma.
\begin{lemma}\label{l:l2}
Let $Y = \{Y_n\}$ be a dyadic supermartingale with drift $d>0$,
\begin{equation}\notag
\frac12(Y(\beta^+) + Y(\beta^-)) \leq Y(\beta) - d,
\end{equation}
and its differences are bounded from above
\begin{equation}\notag
Y_n-Y_{n-1} \leq C.
\end{equation}
Then for any $k\geq0$ and $r\leq d/10$ one has
\begin{equation}\notag
\sharp\{\beta\in E: \ d(\beta,\omega) = k,\; Y(\beta) \geq -kr\} \leq C(r)2^{\frac{k}{2}}.
\end{equation}
\end{lemma}
Assume for a moment that we have Lemma \ref{l:l2}, and hence \eqref{e:71}. Then, we get
\begin{equation}\label{e:64}
.
\end{equation}
\begin{equation}\notag
\begin{split}
&\sum_{\alpha\in E}\varphi^{\delta}(\alpha)\sum_{\beta\in B(\alpha)\bigcap A(\alpha)}\varphi(\beta)|\beta|= \sum_{\alpha\in E}\varphi^{\delta}(\alpha)\sum_{k\geq0}\sum_{\beta\in B(\alpha,k)\bigcap A(\alpha,k)}\varphi(\beta)|\beta|\\
&\leq \sum_{\alpha\in E}\varphi^{\delta}(\alpha)\sum_{k\geq0}\sum_{\beta\in B(\alpha,k)}(1+r)^k\varphi(\alpha)|\beta|=\sum_{\alpha\in E}\varphi^{\delta+1}(\alpha)\sum_{k\geq0}|B(\alpha,k)|(1+r)^k|\beta|\\
&\leq \sum_{\alpha\in E}\varphi^{\delta+1}(\alpha)\sum_{k\geq0}|B(\alpha,k)|2^{2rk}|\alpha|2^{-k}\leq\sum_{\alpha\in E}\varphi^{\delta+1}(\alpha)|\alpha|\sum_{k\geq0}|B(\alpha,k)|2^{-k(1-2r)}\\
&\lesssim \sum_{\alpha\in E}\varphi^{\delta+1}(\alpha)|\alpha|\sum_{k\geq0}2^{\frac{k}{2}}2^{-k(1-2r)}\lesssim\sum_{\alpha\in E}\varphi^{\delta+1}(\alpha)|\alpha|,
\end{split}
\end{equation}
and we are done.
\end{proof}

Lemma \ref{l:l2} clearly follows from the following rescaled driftless version.
\begin{lemma}\label{l:l3}
Let $X = \{X_n\}$ be a dyadic supermartingale,
\begin{equation}\label{e:80}
\frac12(X(\beta^+) + X(\beta^-)) \leq X(\beta),
\end{equation}
and its differences are bounded from above
\begin{equation}\label{e:81}
X_n-X_{n-1} \leq 1.
\end{equation}
Then for any $k\geq0$ and $\eta>0$ one has
\begin{equation}\label{e:82}
\sharp\{\beta\in E: d(\beta,\omega) = k,\; X(\beta) \geq k\eta\} \leq C(\eta)2^{\frac{k}{2}}.
\end{equation}
\end{lemma}
This Lemma is in turn a corollary of Azuma-Hoeffding inequality (essentially a good-$\lambda$ argument for supermartingales).
\begin{proposition}[Azuma-Hoeffding inequality]
Let $Z = \{Z_n\}$ be a supermartingale with bounded differences,
\[
|Z_n-Z_{n-1}| \leq c_n,\quad n\in\mathbb{N}.
\]
Then
\begin{equation}\label{e:83}
\mathbb{P}(Z_n-Z_0 \geq K) \leq e^{\frac{-K^2}{2\sum_{k=1}^nc_k^2}},\quad  K>0,\; n\in\mathbb{N}.
\end{equation}
\end{proposition}
While the proposition above requires the supermartingale differences to be bounded above and below, it is not really  relevant here. Namely, assume $X$ satisfies the hypothesis of Lemma \ref{l:l3} and let $S = \{S_n\}$ be its differences,
\[
S(\beta) = X(\beta) - X(P(\beta)),\quad \beta\in E,
\] 
where $P(\beta)$ is the parent of $\beta$. By \eqref{e:81}, $S \leq 1$. Consider the set 
\[
F = \{\beta\in E: S(\beta)\leq -2\},
\]
and define
\begin{equation}\notag
\tilde{S}(\beta) =\begin{cases}
-2,\quad &\beta\in F,\\
S(\beta)\quad &\textup{otherwise}.
\end{cases}
\end{equation}
Now let 
\[
Z(\beta) = \sum_{\alpha \supseteq \beta}\tilde{S}(\alpha).
\]
Clearly $Z$ is still a dyadic supermartingale, since $\tilde{S}(\beta^-) + \tilde{S}(\beta^+) \leq 0 $ by \eqref{e:80}. Also $\tilde{S} \geq S$, hence $Z \geq X$. It is easily seen that that $Z$ has bounded differences, and therefore satisfies Azuma-Hoeffding inequality.

\section{Conformally invariant Hardy's inequality}\label{sec: conf inv}

While the right hand side of the Hardy's inequality \eqref{HardyInequality} does not depend upon the choice of the root vertex $o$,
the Hardy operator contained in the left hand side does, and consequently also the optimal constant $[\mu]=[\mu]_o$ depends on this choice. It is therefore natural to seek an alternative ``conformal'' invariant theory. The term ``conformal invariant'' should be interpreted in the sense that as \eqref{HardyInequality} corresponds, as explained in the introduction, to a Carleson inequality for Besov spaces, in the same way the inequality we are going to introduce should correspond to a continuous inequality which remains invariant under the group of automorphisms of the unit disc.

We consider here the case $p = 2$, $\pi \equiv 1$. We also assume the tree is dyadic and not
rooted: each vertex is the endpoint of three edges, and $T$ is endowed with a
rich group of automorphisms which, having the
Poincar{\'e} distance in mind, play in $T$ the role of conformal
automorphisms. Such automorphisms are also isometries with respect to the distance $d$ and act naturally also on the boundary $\partial T$ (see \cite{figa-talamanca} for a comprehensive exposition on the topic). Once we fix a root $o$, there are $3 \times 2^{n - 1}$ vertices
at distance $n$ from it.

It is easily seen that the Hardy's inequality \eqref{HardyInequality}, holding for functions $f:E\to \mathbb{R}$,
is equivalent to
\begin{equation}\label{eq:HardyVertex}
 \int_{\overline{T}} | F (x) |^2 \mu (x) \leq [\mu]_o \left( | F (o) |^2 +
   \sum_{\alpha\in E} | \nabla F (\alpha) |^2 \right) , \quad F:\overline{T}\to \mathbb{R},   
\end{equation}
where $\nabla F (\alpha) = F (e (\alpha)) - F (b (\alpha))$ depends on the
choice of the root, but $| \nabla F (\alpha) |^2$ does not. A first attempt to
write down a ``conformally invariant'' formulation of the Hardy's inequality is,
assuming that $\mu (T) = 1$,
\begin{equation}\tag{CH}\label{eq:ConformalHardy}
  \int_{\overline{T}} | F (x) - \mu (F) |^2 d \mu (x) \leq [\mu]_{inv} \sum_{\alpha} |
   \nabla F (\alpha) |^2,   
\end{equation}
where $\mu (F) = \int_T F d \mu$ is the mean of $F$ and $[\mu]_{inv} \in [0, + \infty]$ the best
constant in the inequality. The invariance is the
following. 

Let $\Psi$ be an isometry of $\overline{T}$ and define $\Psi_{\ast} \mu (A) =
\mu (\Psi^{-1} (A))$ and $\Psi^{\ast} F (x) = F (\Psi (x))$, $A\subseteq \overline{T}$. Then,
\begin{eqnarray*}
  \int_{\overline{T}} | F (x) - \Psi_{\ast} \mu (F) |^2 d \Psi_{\ast} \mu (x) & = & \int_{\overline{T}}
  | \Psi^{\ast} F (y) - \mu (\Psi^{\ast} F) |^2 d \mu (y)\\
  & \leq & [\mu]_{inv} \sum_{\beta} | \nabla \Psi^{\ast} F (\beta) |^2 =  [\mu]_{inv} \sum_{\alpha} | \nabla F (\alpha) |^2,
\end{eqnarray*}
showing that $[\Psi_{\ast} \mu]_{inv} = [\mu]_{inv}$.

Observe that the finiteness of $[\mu]_o$ in \eqref{eq:HardyVertex} implies that $\mu$ has no atoms on $\partial T$. On the other hand, if $\mu$ is a Dirac delta measure supported on the boundary the left hand side of \eqref{eq:ConformalHardy} vanishes, while the average $\mu
= \frac{\delta_x + \delta_y}{2}$ of two Dirac delta gives a true, non trivial inequality.

We will show that if $\mu$ is not a boundary Dirac delta, then \eqref{eq:ConformalHardy} is equivalent to \eqref{HardyInequality}. We present two separate arguments, one for measures supported on the boundary of the tree, and another one for measures supported on the vertex set. The first case is proved by means of the isocapacitary characterization. Since also the capacity $\tmop{Cap}(A)=\tmop{Cap}_{2,\pi}(A)$ of a set $A\subseteq\overline{T}$ depends on the choice of the root $o$, in this section we will denote it by $\tmop{Cap}_o(A)$, making explicit the dependence so far kept implicit. On the other hand, \textit{condensers capacity} is invariant under Möbius trasformations of the unit disc and, on trees, under the action of automorphisms. Given two disjoint sets $A,B\subseteq\partial T$, each being a finite union of arcs\footnote{A tent in a non-rooted dyadic tree is any rooted dyadic subtree, and an arc is the boundary of a tent. We
consider, here and in the whole section, finite unions of tents only in order to avoid the complication
of properties which, in the general case, only hold outside sets of capacity
zero.}, we define the capacity of the condenser $(A,B)$ as
\begin{equation*}
    \tmop{Cap} (A, B) = \inf \left\{ \sum_{\alpha \in E } |\nabla F
     (\alpha)|^2 : F |_A = 1, \ F |_B = 0 , \text{ and } |
     \nabla F | \text{ has finite support} \right\}.
\end{equation*}
The next result shows that, for measures supported on the boundary, \eqref{eq:HardyVertex} and \eqref{eq:ConformalHardy} are equivalent, and relates the optimal constants to a capacitary expression.

\begin{theorem}\label{thm:hardy equiv hardy inv}
  Let $\mu \geq 0$ be a Borel probability measure on $\overline{T}$, giving no mass to vertices, and not
  being a Dirac delta on the boundary. Then,
  \[ [[\mu]] \approx \inf_{o \in T} [\mu]_o \approx \sup_{A, B \subseteq \partial T}
     \frac{\mu (A) \mu (B)}{\tmop{Cap} (A, B)}, \]
  where the supremum is over all couples of finite union of arcs.
\end{theorem}

The following lemma provides a recursive formula for calculating the capacity of a condenser of a special type. This kind of formulas arise often in the setting of discrete capacities.

\begin{lemma}\label{lem: capacity formulas}
  Let $o \in T$ be a vertex, and let $T_1, T_2, T_3$ be the dyadic subtrees having it
  as pre-root. Let $A_j \subseteq \partial T_j$ be finite union of arcs. Then,
  \[ \tmop{Cap} (A_1,  A_2 \cup A_3) = \frac{\tmop{Cap}_o (A_1)\cdot(\tmop{Cap}_o (A_2) + \tmop{Cap}_o
     (A_3))}{\tmop{Cap}_o (A_1) + \tmop{Cap}_o (A_2)
     + \tmop{Cap}_o (A_3)} . \]
In particular, for $i\neq j\neq k$, $ \tmop{Cap} (A_i,  \partial T_j \cup \partial T_k)\approx \tmop{Cap}_o (A_i)$.
\end{lemma}

\begin{proof}
Let $c_j=\tmop{Cap}_o (A_j)$.  As in \cite[Proposition 1]{Arcozzi2016}, it can be proved that there exists an extremal
  function $F \geq 0$ on $T$ such that (a) $\lim_{T \ni x \rightarrow
  \zeta} F (x) = 0$ for $\zeta \in A_2 \cup A_3$; (b) \ $\lim_{T \ni x
  \rightarrow \xi} F (x) = 1$ for $\xi \in A_1$; (c) $\sum_{\alpha \in E}
  \nabla F (\alpha)^2 = \tmop{Cap} (A_1,  A_2 \cup A_3)$. Of course, such $|
  \nabla F |$ is not finitely supported. Similarly, there exist analogous
  extremal functions $F_j$ for $c_j$, $j = 1, 2, 3,$ with
  $F_j (a) = 0$; \ $\lim_{T \ni x \rightarrow \xi} F_j (x) = 1$ for $\xi
  \in A_j$; $\sum_{\alpha \in E} \nabla F_j (\alpha)^2 = c_j$.
  
  By estremality, it is obvious that there are numbers $0<s_j\leq 1$ such that $| \nabla F (\alpha) | = s_j | \nabla
  F_j (\alpha) |$ on the edges $\alpha$ of $T_j$, and since $| \nabla F |$
  adds to one along geodesics going from $A_1 \cup A_2$ to $A_3$, it must be $s_1 + s_2 = s_1
  + s_3 = 1$. Again by minimality, we look for $t = s_1$ which minimizes
  \begin{eqnarray*}
    f (t) & = & \| \nabla F \|_{\ell^2}^2 = (1 - t)^2 (\| \nabla F_2
    \|_{\ell^2}^2 + \| \nabla F_3 \|_{\ell^2}^2) + t^2 \| \nabla F_1
    \|_{\ell^2}^2\\
    & = & (1 - t)^2 (c_2 + c_3) + t^2
    c_1,
  \end{eqnarray*}
  which is minimal when $t = \frac{c_2 + c_3}{c_1 + c_2 + c_3}$, thus
  \[ \tmop{Cap} (A_1,  A_2 \cup A_3) = f \left( \frac{c_2 +
     c_3}{c_1 + c_2 +
     c_3} \right) = \frac{c_1\cdot(c_2 + c_3)}{c_1 + c_2
     + c_3} . \]

Since clearly $\tmop{Cap}_{o} (T_1)=\tmop{Cap}_{o} (T_2)=\tmop{Cap}_{o} (T_3)=1$ and $\tmop{Cap}_o (A_i)\leq 1$,
\begin{equation*}
    \tmop{Cap} (A_i,  \partial T_j \cup \partial T_k)=\frac{2 \tmop{Cap}_o (A_i)}{2 + \tmop{Cap}_o (A_i)}
  \approx \tmop{Cap}_o (A_i).  
\end{equation*}
\end{proof}

\begin{proof}[Proof of Theorem \ref{thm:hardy equiv hardy inv}]
In one direction,
   \begin{eqnarray*}
    \int_{\partial T} \big(F - \mu (F)\big)^2 d \mu & = & \int_{\partial T} \big((F - F
    (o)) - \mu (F - F (o))\big)^2 d \mu\\
    & = & \int_{\partial T} \big(F - F (o)\big)^2 d \mu - \left( \int_{\partial T} (F
    - F (o)) d \mu \right)^2\\
    & \leq & [\mu]_o \cdot \sum_{\alpha \in E} \nabla F (\alpha)^2,
  \end{eqnarray*}
  hence, $ [[\mu]] \leq \inf_{o \in T} [\mu]_o$.
  
 In the other direction, consider closed subsets $A, B \subseteq \partial T$,
  which we might assume to be finite unions of arcs, and a function $F$ with finitely supported $|\nabla F|$ such that $F = 1$ on $A$ and $F = 0$ on $B$. Then,
  \begin{eqnarray*}
    {}[[\mu]] \cdot \sum_{\alpha \in E} \nabla F (\alpha)^2 & \geq &
    \int_{\partial T} | F - \mu (F) |^2 d \mu\\
    & = & \frac{1}{2}\int_{\partial T\times\partial T}  \big(F (x) - F
    (y)\big)^2d \mu (x)d \mu (y)\\
    & \geq & \mu (A) \mu (B) .
  \end{eqnarray*}
  Passing to the infimum over all such $F$'s, we obtain
  \[ \mu (A) \mu (B) \leq [[\mu]] \cdot \tmop{Cap} (A, B), \]
  giving $[[\mu]] \geq \sup_{A, B \subseteq \partial T} \frac{\mu (A) \mu
  (B)}{\tmop{Cap} (A, B)}$.
  
  It is not difficult to see that finiteness of $\sup_{A, B \subseteq
  \partial T} \frac{\mu (A) \mu (B)}{\tmop{Cap} (A, B)}$ implies that $\mu$ does not
  have more than an atom on the boundary, which would then be a Dirac delta. Hence, if $\mu$ has boundary atoms the statement holds.
  
  Suppose now that $\mu$ is atomless.  We claim:
  \begin{itemize}
      \item[(i)] if $\mu$ is a
  probability measure on $\partial T$ having no atoms, then there are disjoint arcs $I_1 \cup I_2 = \partial T$, such that $1 / 3 \leq \mu (I_j) \leq 2 / 3$.
   \end{itemize}
  
  With the claim given, let $A$ be a finite union of arcs, let $I_1, I_2$ as given by (i), and $A_1=A\cap I_2$, $A_2=A\cap I_1$. Let $o_j$ be the pre-root of the dyadic tree $T_1^j$ having boundary $I_j$ and $T_2^j,T_3^j$ the two other dyadic trees having it as a pre-root, so that, for $k\neq j$, $I_k=\partial T_2^j \cup \partial T_3^j$, and,
 \begin{equation*}
     [[\mu]]\geq \sup_{A, B \subseteq \partial T} \frac{\mu (A) \mu (B)}{\tmop{Cap} (A, B)} \geq \frac{\mu (A_j) \mu
     (I_j)}{\tmop{Cap} (A_j, I_j)} \approx \frac{\mu (A_j)}{\tmop{Cap}_{o_k} (A_j)},
 \end{equation*}
  Observing that $\tmop{Cap}_{o_k} (A_j)\leq \tmop{Cap}_{o_j} (A_j)$, for $o\in \lbrace o_1,o_2\rbrace$ we have
  \begin{equation*}
      \frac{\mu (A)}{\tmop{Cap}_o(A)}\leq \frac{\mu (A_1)}{\tmop{Cap}_{o_2} (A_1)}+\frac{\mu (A_2)}{\tmop{Cap}_{o_1} (A_2)}\leq 2 \sup_{A, B \subseteq \partial T} \frac{\mu (A) \mu (B)}{\tmop{Cap} (A, B)}\leq 2  [[\mu]].
  \end{equation*}
  
  By the isocapacitary condition \eqref{eq:isocapacitary} characterizing measures satisfying the Hardy's
  inequality, $[\mu]_o=\sup\lbrace \frac{\mu (A_j)}{\tmop{Cap}_o (A_j)}: \  A \text{ finite union of arcs} \rbrace$, and we have the promised statement.
  
  We now come to the proof of the claim. Choose a vertex $x_0\in T$ and for $j=1,2,3$  let $T_j^0$ be the dyadic subtrees having pre-root at it and call $x_j$ the neighborhood point of $x_0$ lying in $T_j^0$. For at least one $j$, $\mu (\partial T_j^0) \geq 1 / 3$; say $j =
  1$. If $\mu (\partial T_1^0) \leq 2 / 3$ we set $I_1 = \partial T_1^0$
  and we are done. Otherwise, let $T_1^1, T_2^1$ be the two infinite subtrees of $T_1^0$ with pre-root in $x_1$. For one $j\in \lbrace 1,2\rbrace$ we have $\mu (\partial T_j^1) \geq 1 / 3$;
  let it again be $j = 1$. As before, set $I_1 = \partial T_1^1$ if $\mu
  (\partial T_1^1) \leq 2 / 3$, otherwise consider the two infinite subtrees of $T_1^1$ rooted at the neighboroods of $x_1$, and iterate the reasoning. If there is no stop,
  we have a family of nested tents $\partial T_1^0 \supset \partial T_1^1
  \supset \ldots$, whose intersection is a single boundary point $x$ with
  $\mu (\{ x \}) \geq 2 / 3$, contradicting the assumption.
  
\end{proof}

We come now to the case of measures supported on the vertex set of $T$. Since no extra work is required, we present a proof which holds in much higher generality.

\begin{proposition}
Let $X$ be a locally compact space and $B$ be a Banach space of functions on
$X$ such that for every compact $K$ there is $C (K)$ with $\sup_{x \in K} | f
(x) | \leq C (K)$ if $\| f \|_B \leq 1$. Then the following are
equivalent for a probability measure $\mu$:
\begin{itemize}
  \item[(a)] $\int | f - \mu (f) |^p d \mu \leq C_0 \| f \|_B^p$,
  
  \item[(b)] $\int | f |^p d \mu \leq C_1 \| f \|_B^p .$
\end{itemize}
\end{proposition}

\begin{proof}

If (b) holds, then
\begin{equation*}
  \| f - \mu (f) \|_{L^p (\mu)}  \leq  \| f \|_{L^p (\mu)} + | \mu (f)
  |\leq  2 C_1 \| f \|_B .
\end{equation*}
Suppose (a) holds and (b) does not. By (a):
\[ \| f \|_{L^p (\mu)} \leq \| f - \mu (f) \|_{L^p (\mu)} + | \mu (f) |
   \leq C_0 \| f \|_B + | \mu (f) |, \]
then there exists a sequence $f_n$ in $B$ with (i) $\| f_n \|_B = 1$, (ii) $\|
f_n \|_{L^p (\mu)} \nearrow \infty$, (iii) $M_n \assign | \mu (f_n) | \geq 2^{-1}
\| f_n \|_{L^p (\mu)}$.

We can find compact sets $K_n \nearrow X$ such that $\int_{K_n} | f_n | d \mu
> \frac{M_n}{2}$. For any fixed compact $S$ we have that $\int_S
| f_n | d \mu\leq C (S) \leq \frac{M_n}{4}$ for $n \geq n (S)$, hence,
\[ \int_{K_n / S} | f_n | d \mu \geq \frac{M_n}{2} - \frac{M_n}{4} =
   \frac{M_n}{4} \]
for $n \geq n (S)$. Thus
\begin{equation*}
  \frac{M_n^p}{4^p}  \leq  \left( \int_{K_n / S} | f_n | d \mu
  \right)^p \leq  \| f \|_{L^p (\mu)}^p \mu (X \setminus S)^{p / p^*},
\end{equation*}
i.e., $\frac{| \mu (f_n) |}{\| f \|_{L^p (\mu)}} \leq 4 \mu (X \setminus
S)^{1 / p^*}$, which can be made as small as we wish, leading to a contradiction. 
\end{proof}

\begin{problem}\label{prob:conf inv strong}
A more general problem is that of characterizing the probability measures
$\lambda$ on $T \times T$ such that
\begin{equation}\tag{a}\label{(a)}
   \quad \int_T \int_T | F (x) - F (y) |^2 d \lambda (x, y) \leq \{
   \lambda \} \sum_{\alpha} | \nabla F (\alpha) |^2, 
\end{equation}
which is as well conformally invariant, if we set $\Psi_{\ast} \lambda (A
\times B) = \lambda (\Psi (a) \times \Psi (B))$, and reduces to the above for
$\lambda = \mu \otimes \mu$:
\begin{eqnarray*}
  \int_T \left| F (x) - \int_T F (y) d \mu (y) \right|^2 d \mu (x) & = &
  \int_T \left| \int_T (F (x) - F (y)) d \mu (y) \right|^2 d \mu (x)\\
  & \leq & \int_T \int_T | F (x) - F (y) |^2 d \mu (y) d \mu (x)\\
  & = & 2 \left[ \int_T F (x)^2 d \mu (x) - \left( \int_T F (x) d \mu (x)
  \right)^2 \right]\\
  & = & 2 \left[ \int_T \left| F (x) - \int_T F (y) d \mu (y) \right|^2 d \mu
  (x) \right] .
\end{eqnarray*}

No characterization of the measures $\lambda$ for which (a) holds is known,
to the best of our knowledge.
\end{problem}
\begin{problem}
A natural generalization of the conformal Hardy's inequality is to consider a $p$-version, or even a weighted version of it. Namely, given a weight function $\pi : E \to \mathbb{R}_+$ characterize positive Borel probability measures such that there exists a constant $[\mu]_{inv} $ possibly depending on $\mu, \pi,  $ and $\pi$ such that 
\begin{equation}
  \int_{\overline{T}} | F (x) - \mu (F) |^p d \mu (x) \leq [\mu]_{inv} \sum_{\alpha} |
   \nabla F (\alpha) |^p \pi(\alpha),   
\end{equation}
for all $F: \overline{T} \to \mathbb{R}. $
\end{problem}
There is a related interesting, conformally
invariant interpolation problem.

\begin{problem}
Set $\| F \|^2 = \sum_{\alpha \in E}
\nabla F (\alpha)^2$. Then, $| F (x) - F (y) | \leq \| F \| \cdot d (x,
y)^{1 / 2}$. Given a subset $Z \subset T$, we say that it is
{\tmem{universally interpolating}} for the seminorm $\| \cdot \|$ if
\begin{enumerateroman}
  \item $\sum_{z \in T} \sum_{w \in T} \frac{| F (z) - F (w) |^2}{d (z, w)}
  \leq C \| F \|^2,$
  
  \item for all sequences $\{ a (z) \}_{z \in Z}$ such that $\sum_{z \in T}
  \sum_{w \in T} \frac{| a (z) - a (w) |^2}{d (z, w)} < \infty$ there exists
  $F$ with $\| F \| < \infty$ such that $F (z) = a (z), \forall z \in Z$.
\end{enumerateroman}
Condition (i) says that the measure $\lambda = \sum_{z, w \in T}
\frac{\delta_{(z, w)}}{d (z, w)}$ satisfies (a). We call the sequence 
{\tmem{onto interpolating}} if just (ii) holds. We think that it is an
interesting problem finding characterizations of universally, or onto,
interpolating sequences.
\end{problem}

\section{Miscellaneous results}\label{sec:miscellaneous}
\subsection{Compactness} 

We will briefly discuss the compactness conditions for the Hardy operator. As it is natural to expect the compactness of the Hardy operator corresponds to the ``vanishing" versions of the conditions characterizing boundedness.

\begin{theorem}\label{comppp}
 The map $\III:\ell^p(\pi)\to L^p(\mu)$ is compact if and only if
 \begin{equation}\label{eqcomppp}
\lim_{n\to\infty}\sup_{d(\alpha)\geq n}\frac{ \sum_{\beta\subseteq \alpha }\pi(\beta)^{1-p^*}\mu(S(\beta))^{p^*}}{\mu(S(\alpha))}=0.  
 \end{equation}
\end{theorem}
\begin{proof}
 Suppose \eqref{eqcomppp} holds, fix $n\ge1$ and consider the finite rank truncation of $\III^\ast_\mu$,  
 $\III^\ast_{\mu,n}g(\alpha)=\chi(
 d(\alpha) \leq n)\III^\ast_{\mu}g(\alpha)$. Then,
 \begin{eqnarray*}
  \left\|\III^\ast_{\mu}g-\III^\ast_{\mu,n}g\right\|_{\ell^{p^*}(\pi^{1-p^*})}^{p^*}&=&
  \sum_{d(\alpha)=n}\sum_{\beta\subseteq \alpha}\left|\III^\ast_{\mu}g(\beta)\right|^{p^*}\pi(\beta)^{1-p^*}\crcr
  &\le &C\sum_{d(\alpha)=n}\sup_{\gamma\subseteq \alpha}\frac{ \sum_{\beta\subseteq \gamma}\pi^{1-p^*}(\beta)\mu(S(\beta))^{p^*}}{\mu(S(\gamma))}\|g\chi_{S(\alpha)}\|_{\ell^{p^*}(\mu)}^{p^*}\crcr
  &\le&C\sup_{d(\alpha)\geq n}\frac{ \sum_{\beta\subseteq \alpha}\pi^{1-p^*}(\beta)\mu(S(\beta))^{p^*}}{\mu(S(\alpha))}\|g\|_{\ell^{p^*}(\mu)}^{p^*}\to0.
 \end{eqnarray*}
 Therefore $\III^*_\mu$ is compact as an operator norm limit of finite rank operators. By Schauder's Theorem \cite[p. 243 Theorem 7]{lax2002functional} $\III$ is also compact. 
 
 To see the necessity of the condition we can again work with $\III^*_\mu$ courtesy of Schauder's Theorem. Suppose that \eqref{eqcomppp} does not hold. Then there exists an $\varepsilon>0$ and a sequence of edges $\{ \alpha_k \}_k$ such that $\lim_k d(\alpha_k)= \infty$ such that 
 \[ \frac{ \sum_{\beta\subseteq \alpha_k}\pi^{1-p^*}(\beta)\mu(S(\beta))^{p^*}}{\mu(S(\alpha_k))} \geq \varepsilon. \]
 
Then consider the sequence of testing functions $g_k:=\mu(S(\alpha_k))^{-1/p^*}\chi_{S(\alpha_k)}$, which converges weakly to $0$ in the space $L^{p^*}(\mu).$ By compactness we must have 
\[ 0 = \lim_k \norm{\III^*_\mu g_k}_{\ell^{p^*}(\pi^{1-p^*})}^{p^*} = \frac{ \sum_{\beta\subseteq \alpha_k}\pi^{1-p^*}(\beta)\mu(S(\beta))^{p^*}}{\mu(S(\alpha_k))},  \]
which contradicts that the above quantity is bounded below by $\varepsilon.$
\end{proof}

In a very similar way one can characterize the compactness of the Hardy operator in terms of a {\it vanishing capacitary condition}. We state the theorem without a proof in order to avoid repetition. 

\begin{theorem}\label{comp-cap}
  The map $\III:\ell^p \to L^p(\mu)$ is compact if and only if 
  \[ \lim_n \sup_{d(\alpha_i)\geq n} \Capac_{p,\pi}\Big( \bigcup_{i=1}^k S(\alpha_i) \Big)^{-1} \sum_{i=1}^k\mu(S(\alpha_i)) = 0. \]
  
\end{theorem}

\begin{problem}
Find a  simple characterization of measures $\mu$ such that the Hardy operator $\III:\ell^2(\pi) \to L^2(\mu) $ belongs to the $p-$Schatten ideal. 

In \cite{Luecking1987} Luecking has studied trace ideal criteria for Toeplitz operators on weighted Bergman spaces. It is very possible that some of his results apply also to our case, although his results are not complete. 
\end{problem}

\subsection{Not sufficiency of the simple box condition} \label{sec:simplebox}

A necessary condition for the Hardy's inequality to hold is the simple, one-box condition
\begin{equation} 
 \sup_{\alpha \in E}\mu(S(\alpha)) d_\pi(e(\alpha))^{p-1} := [[\mu]]_{sc} < +
\infty,
\end{equation}
where $d_{\pi} (x) = \sum_{\alpha \in [o^*,x]} \pi (\alpha)^{1 - p^*}$. In
fact, if $\varphi = \chi_{[o^*,x]} \pi^{1 - p^*}$, then
$\int_{\overline{T}}\III\varphi^p d\mu\geq\mu((S(\alpha))d_{\pi} (e(\alpha))^p$, 
while $\| \varphi \|_{L^p
(\pi)}^p = d_{\pi} (e(\alpha))$.

In general, however, \eqref{CarlesonBox} is not sufficient. We provide here a counterexample for $p =
2$, $\pi = 1$, and $T$ the dyadic tree. Working a bit, the example can be
modified to hold for $1 < p < \infty$.

\begin{example}[\eqref{CarlesonBox} $\centernot \implies$ \eqref{HardyInequality}]\label{counterexample}
 Let $T$ be a dyadic tree, $\pi\equiv 1$ and $p=2$. At level $N_k = 2^k k$, choose $2^k$ vertices $\lbrace z_n^k\rbrace_{n=1}^{2^k}$, with children named $x_n^k$ and $y_n^k$, such that $z_{2 n - 1}^k \wedge z_{2n}^k=y_n^{k - 1}$. Choose points $w_n^k\in S(x_n^k)$ with $d (w_n^k) = M_k =2^k k^2$.

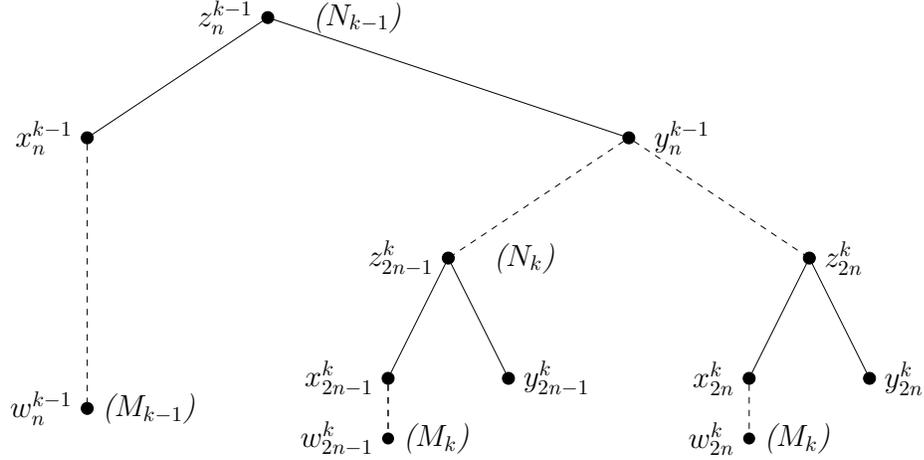
\begin{figure}[H]
    
\begin{center}

\scalebox{0.8}{%
\tikzstyle{every node}=[font=\large]
\begin{tikzpicture}

\node[shape=circle,draw=black, fill=black, scale=0.5, label=left:{$z_n^{k-1}$}, label={right:{\quad ($N_{k-1}$)}}](A) at (0,0) {};

\node[shape=circle,draw=black, fill=black, scale=0.5, label=left:{$x_n^{k-1}$}](Aa) at (-3,-2) {};
\node[shape=circle,draw=black, fill=black, scale=0.5, label=right:{\ $y_n^{k-1}$}](Ac) at (6,-2) {};

\node[shape=circle,draw=black, fill=black, scale=0.5, label=left:{$w_n^{k-1}$}, label={right:{ ($M_{k-1}$)}}](Aaa) at (-3,-6.5) {};

\node[shape=circle,draw=black, fill=black, scale=0.5, label=left:{$z_{2n-1}^k$}, label={right:{\quad ($N_{k}$)}}](Aca) at (3,-4) {};

\node[shape=circle,draw=black, fill=black, scale=0.5, label=right:{$z_{2n}^k$}](Acb) at (9,-4) {};


\node[shape=circle,draw=black, fill=black, scale=0.5, label=left:{$x_{2n-1}^k$}](Acaa) at (2,-6) {};
\node[shape=circle,draw=black, fill=black, scale=0.5, label=right:{$y_{2n-1}^k$}](Acab) at (4,-6) {};

\node[shape=circle,draw=black, fill=black, scale=0.5, label=left:{$x_{2n}^k$}](Acba) at (8,-6) {};
\node[shape=circle,draw=black, fill=black, scale=0.5, label=right:{$y_{2n}^k$}](Acbb) at (10,-6) {};


\node[shape=circle,draw=black,draw=none, fill=black, scale=0.5, label=left:{$w_{2n-1}^k$}, label={right:{ ($M_k$)}}](Acaaa) at (2,-7) {};
\node[shape=circle,draw=black,draw=none, fill=black, scale=0.5, label=left:{$w_{2n}^k$}, label={right:{ ($M_k$)}}](Acbaa) at (8,-7) {};

  
  \path [- ] (A) edge node[left] {} (Aa); 
  \path [-] (A) edge node[right] {} (Ac); 

 \path [dashed] (Aa) edge node[left] {} (Aaa); 

 \path [dashed] (Ac) edge node[left] {} (Aca); 
 \path [dashed] (Ac) edge node[right] {} (Acb); 
 
 \path [- ] (Aca) edge node[left] {} (Acaa);  \path [- ] (Aca) edge node[left] {} (Acab);
 \path [- ] (Acb) edge node[left] {} (Acba);  \path [- ] (Acb) edge node[right] {} (Acbb);
 
 \path [dashed] (Acaa) edge node[left] {} (Acaaa);
 \path [dashed] (Acaa) edge node[left] {} (Acaaa);  
 
  \path [dashed] (Acba) edge node[left] {} (Acbaa);  

\end{tikzpicture}
}

\end{center}
\caption{A snapshot from the defined family of points and their relations: continuous lines represent edges and dashed lines paths of (many) edges. Next to some vertices is specified, in parenthesis, their distance from the origin.\footnotesize}
    \label{fig:my_label}
\end{figure}

Let $\mu = \sum_{k, n} \frac{1}{M_k} \delta_{w_n^k}$. A simple reasoning shows that it suffices to verify the one-box condition at the nodes $z_n^k$:

\begin{eqnarray*}
  \mu (S (z_n^k)) & = & \sum_{j = 0}^{\infty} \frac{2^j}{M_{j + k}} = 2^{- k}
  \sum_{j = 0}^{\infty} \frac{1}{(j + k)^2} \approx \frac{1}{N_k}=d (z_n^k)^{-1}.
\end{eqnarray*}
On the other hand, the mass--energy condition \eqref{EqMassEnergy} fails. To see this, denote by $Z$ the minimal subtree containing all the $z_n^k$ points. Then,
\begin{equation*}
\begin{split}
    \sum_{x\in T}\mu(S(x))^2 &\geq \sum_{x\in Z}\mu(S(x))^2=\sum_{k=0}^\infty(N_{k+1}-N_k)2^k\Big(\sum_{j = 0}^{\infty} \frac{1}{(j + k)^2}\Big)^2\\
    &\approx \sum_{k=0}^\infty\Big(\frac{N_{k+1}}{N_k^2}-\frac{1}{N_k}\Big)2^k\approx \sum_{k=0}^\infty\frac{1}{k}=+\infty.
\end{split}
\end{equation*}
\end{example}

\begin{problem}
Find  a characterization of those couples $(p, \pi)$ for which
\eqref{CarlesonBox} is not sufficient on the dyadic tree.
\end{problem}


\subsection{Two opposite examples}

In the generality in which we have stated it, the dyadic Hardy's inequality
covers a variety of contexts; some of them very rich, some very poor. The
richest context is in our opinion the unweighted case $\pi = 1$, corresponding
to the classical Dirichlet space. We consider here two cases at the opposite
extremes.

\begin{example}[Boundary having null capacity]
Consider the dyadic, infinite tree
with the weight $\pi (\alpha) = 2^{- d (\alpha)}$ and $p = 2$. The reader
familiar with martingale theory can fruitfully think of $\pi$ as the
probability of a fair coin tossed $d (\alpha)$ times. Let us show that with
this choice $\tmop{Cap}_{2,\pi} (\partial T) = 0$. Let $g_N (\alpha) = \frac{1}{N}$ if $d (\alpha) \le N$
  and $g_N (\alpha) = 0$ elsewhere. Then,
\begin{eqnarray*}
  \tmop{Cap}_{2,\pi} (\partial T) & = & \inf \left\{ \| f \|^2_{\ell^2 (\pi)} : \ f \ge 0, \ \mathcal{I}f \ge 1 \text{ on } \partial T
  \right\}\\
  & \le & \| g_N \|^2_{\ell^2 (\pi)}= \sum_{n = 0}^N 2^n \frac{1}{N^2} 2^{- n}= \frac{1}{N} \longrightarrow 0, \quad \text{as} \ N\to +\infty.
\end{eqnarray*}
It follows from the isocapacitary inequality \eqref{eq:isocapacitary} that $\partial T$ does not support any Carleson measure $p = 2$ and the chosen weight $\pi$. In particular, the Lebesgue measure $\mu_0 (\partial S (\alpha)) = 2^{- d
(\alpha)}$ does not define a Carleson measure. This fact
is best appreciated having in mind basic martingale theory.

Consider the filtration associated with the infinite tossing of a fair coin,
where $\partial T$ is the probability space and $\mu_0 (\partial S (\alpha)) =
2^{- d (\alpha)}$ is the probability measure. \ A martingale for the
filtration has the form $X_n (\zeta) =\mathcal{I}f (\alpha)$ for $\zeta \in
\partial S (\alpha)$ and $d (\alpha) = n$. The martingale Hardy space
$\mathcal{M}^2$ contains those martingales for which
\begin{eqnarray*}
  \| X \|_{\mathcal{M}^2}^2 & = & \sup_n \int_{\partial T} X_n^2 d \mu_0\\
  & = & \int_{\partial T} \mathcal{I}f^2 d \mu_0\\
  & = & \sum_{\alpha} f (\alpha)^2 \pi (\alpha) < \infty .
\end{eqnarray*}
i.e., $\mathcal{M}^2 = \ell^2 (\pi) \cap \mathcal{M}$, where $\mathcal{M}$ is
the space of all martingales. By definition, $\mu_0$ is Carleson measure for
$\mathcal{M}^2$, but it is not for $\ell^2 (\pi)$. The underlying reason is
that in $\mathcal{M}^2$ cancellations play a prominent role, and this much
enlarges the class of the Carleson measures at the boundary.

But suppose we ask a measure $\mu$ to be Carleson for the variation of the
martingales in $\mathcal{M}^2$,
\[ \int_{\partial T} \left( \sum_{\alpha \in P (\zeta)} | f (\alpha) |
   \right)^2 d \mu (\zeta) \le C \sum_{\alpha} f (\alpha)^2 \pi
   (\alpha), \]
where the variation of $X =\mathcal{I}f$ is $V (X) (\zeta) = \sum_{\alpha \in
P (\zeta)} | f (\alpha) |$. It is easy to see that this is the same as asking
$\mu$ to be Carleson for $\ell^2 (\pi)$, hence $\mu = 0$. This is reflected in
the fact that functions in the classical Hardy space can have unbounded
variation a.e. at the boundary of the unit disc \cite{Rudin1955, Bourgain1993}.
\end{example}

Let us note that, in contrast, for the weights $\pi_{\lambda} (\alpha) =
2^{\lambda d (\alpha)}$, $0 \le \lambda < 1$, the Carleson measures for
$\ell^2 (\pi_{\lambda}) \cap \mathcal{M}$ and $\ell^2 (\pi_{\lambda})$ are the
same \cite{AR2004}. This is in much the spirit of
Beurling's result on exceptional sets \cite{beurling1940ensembles}.

\begin{example}[All boundary points have positive capacity]
It is easy to see that for any tree $T$ and any $\zeta\in \partial T$, $\tmop{Cap}_p^{\sigma} (\{
\zeta \}) = d_{\sigma} (\zeta)^{-1}$. In fact, for any function $f$ which is admissible for $\zeta$,
\begin{equation*}
   \Vert f\Vert_{\ell^p(\sigma)}^p\geq\sum_{\alpha\supset\zeta}|f(\alpha)|^p\sigma(\alpha)\geq\Big(\sum_{\alpha\supset\zeta}\frac{|f(\alpha)|}{d_\sigma(\zeta)^{1/p^*}}\sigma(\alpha)\Big)^p=\frac{(I_\sigma f(\zeta))^p}{d_\sigma(\zeta)^{p-1}}\geq d_\sigma(\zeta)^{1-p},
\end{equation*}
and the right handside is the $\ell^p(\sigma)-$norm of the admissible function taking constant value $d_\sigma(\zeta)^{-1}$ on the edges $\alpha\supset\zeta$ and zero elsewhere. Then, assuming that $\tmop{Cap}_p^{\sigma} (\{
\zeta \}) \ge c > 0$ for all boundary points, is
the same as saying that $\partial T$ is bounded with respect to the distance
$d_{\sigma}$. This is the case, for example, for the weights $\sigma (\alpha) = 2^{\lambda d (\alpha)}$ with
$\lambda > 0$. Under this assumption, all functions $\mathcal{I}f$ with $f \in \ell^p
(\pi)$, are bounded on $\partial T$:
\begin{equation*}
| \mathcal{I}f (\zeta)| \le \left( \sum_{\alpha \in
  [o^*, \zeta]} | f (\alpha)
  |^p \pi (\alpha) \right)^{1 / p} \left( \sum_{\alpha \in [o^*, \zeta]}
  \sigma (\alpha) \right)^{1 / p^*}\le \Vert f\Vert_{\ell^p(\pi)}\Big(\frac{2}{c}\Big)^{1 / p^*},
\end{equation*}
from which follows that  all bounded measures $\mu$ are Carleson for $\ell^p
(\pi)$.

Under this assumption, all functions $\mathcal{I}f$ with $f \in \ell^p
(\pi)$, where as usual $\pi=\sigma^{1-p}$, are continuous up to the boundary, with respect to the natural topology
induced by the visual distance $d(\zeta,\xi)=e^{-d(\zeta\wedge\xi)}$:
\begin{eqnarray*}
  | \mathcal{I}f (\zeta) -\mathcal{I}f (\xi) | & \le & \sum_{\alpha \in
  [\zeta, \xi]} | f (\alpha) |\\
  & \le & \left( \sum_{\alpha \in S (\zeta \wedge \xi)} | f (\alpha)
  |^p \pi (\alpha) \right)^{1 / p} \left( \sum_{\alpha \in [\zeta, \xi]}
  \sigma (\alpha) \right)^{1 / p^*}\\
  & \le & \left( \sum_{\alpha \in S (\zeta \wedge \xi)} | f (\alpha)
  |^p \pi (\alpha) \right)^{1 / p} \Big(\frac{2}{c}\Big)^{1 / p^*},
\end{eqnarray*}
which tends to zero as $\xi \rightarrow \zeta$ by dominated convergence. By
Weierstrass Theorem, all bounded measure $\mu$ are Carleson for $\ell^p
(\pi)$.

\end{example}

\section{Some variations on the structure}\label{sec:variation}

\subsection{The viewpoint of reproducing kernels}\label{sec:rk}
Let us recall to the reader that a \emph{reproducing
kernel Hilbert space (RKHS)} is a Hilbert space $\mathcal{H}$ of functions defined on a set $X$ such that the point evaluation functionals are continuous or, equivalently, such that for any $x\in X$ there exists an element $K_x\in \mathcal{H}$ which fulfills the reproducing property
\begin{equation*}
    f(x)=\langle f, K_x\rangle_{\mathcal{H}}, \quad \text{for all } f\in \mathcal{H}.
\end{equation*}
It is easy to see that the function $K:X\times X\to \mathbb{C}$ given by
\begin{equation*}
    K(y,x):=K_x(y)=\langle K_x, K_y\rangle_{\mathcal{H}},
\end{equation*}
is a kernel on $\mathcal{H}$, that is, it is positive semi-definite, nonzero on the diagonal and satisfies $K(x,y)=\overline{K(y,x)}$. See \cite{aronszajn1950} or \cite{agler2002pick} for a systematic treatment of the topic.

Let $H_K$ be a RKHS of functions on a locally compact space $X$ with continuous reproducing kernel
$K$. A simple and well known ``$T^{\ast} T$ argument'' shows that the imbedding $\tmop{Id} : H_K \rightarrow L^2 (\mu)$ (i.e., $\mu$, a
positive Borel measure on $X$, is a Carleson measure for $H_K$) can be
rephrased various ways in terms of integral inequality on $L^2 (\mu)$.

\begin{lemma}\label{lem:soft}
  Given a RKHS $H_K$ of functions on a locally compact space $X$ with continuous reproducing kernel $K
  (x, y) = K_y (x)$, the following are equivalent for a Borel measure $\mu$ on
  $X$:
  
  \vspace{0.2cm}
  
  \begin{itemize}
    \item[(i)] $\| f \|_{L^2 (\mu)}^2 \leq \{ \mu \}_1 \| f \|_{H_K}^2$, i.e.,
    $\| \tmop{Id}^{\ast}  g \|_{H_K}^2 \leq \{ \mu \}_1 \|  g
    \|_{L^2 (\mu)}^2$,
    
    \vspace{0.2cm}
    
    \item[(ii)] $\langle \tmop{Id}^{\ast}  g, \tmop{Id}^{\ast}  g \rangle_{H_K} =
    \int_X \int_X  g (x) \overline{ g (y)} K_y (x) d \mu (x) d \mu
    (y) \leq \{ \mu \}_1 \|  g \|_{L^2 (\mu)}^2$,
    
      \vspace{0.2cm}
      
    \item[(iii)] $\int_X \int_X  g (x) \overline{ g (y)} \tmop{Re} K_y (x)
    d \mu (x) d \mu (y) \leq \{ \mu \}_2 \|  g \|_{L^2 (\mu)}^2$,
    
      \vspace{0.3cm}
      
    \item[(iv)] $\int_X \int_X  g (x)  g (y) \tmop{Re} K_y (x) d \mu (x) d
    \mu (y) \leq \{ \mu \}_2 \|  g \|_{L^2 (\mu)}^2$, for real (or
    even just positive) $ g$,
    
      \vspace{0.2cm}
      
    \item[(v)] $\| f \|_{L^2 (\mu)}^2 \leq \{ \mu \}_1 \| F \|_{H_{\tmop{Re}
    K}}^2$,
    
      \vspace{0.2cm}
    
    \item[(vi)] $| < \tmop{Id}^{\ast}  g, \psi >_{L^2 (\mu)} | = \left| \int_X
    \int_X  g (x) \overline{\psi (y)} K_y (x) d \mu (x) d \mu (y) \right|
    \leq \{ \mu \}_3 \|  g \|_{L^2 (\mu)} \| \psi \|_{L^2 (\mu)}$,
    
      \vspace{0.2cm}
      
    \item[(vii)]$\| \tmop{Id}^{\ast}  g \|_{L^2 (\mu)}^2 \leq \{ \mu \}_3 \|
     g \|_{L^2 (\mu)}^2$.
  \end{itemize}
  Where the statements are assumed to hold for all $f,g : X \to \mathbb{C}  $ and the constants $\{\mu\}_1, \{\mu\}_2, \{\mu\}_3$, which are the best constants in the respective inequalities, might depend on $\mu$ and $\pi$ but not on $f, g$.
  In particular, $\mu$ is Carleson measure for $H_K$ if and only if it is a
  Carleson measure for $H_{\tmop{Re} K}$.
\end{lemma}

The proof of Lemma \ref{lem:soft} can be found with some variations in different sources \cite[Proposition 4.9]{ARSW2019}, \cite[Lemma 24]{ARS2008}, \cite[pp. 9--10]{ARSW2004}. The proof itself is a short ``soft'' analysis argument, and we write it here for
convenience of the reader. The statement we give is an elaboration of various
statements in the literature.

\begin{proof}[Proof of Lemma \ref{lem:soft}]
(i) says that $\tmop{Id} : H_K \rightarrow L^2 (\mu)$ is bounded with norm
$\{ \mu \}_1^{1 / 2}$, so $\tmop{Id}^{\ast} : L^2 (\mu) \rightarrow H_K$ is
bounded with the same norm, where
\begin{eqnarray*}
  \tmop{Id}^{\ast}  g (x) & = & \langle \tmop{Id}^{\ast}  g, K_x \rangle_{H_K} =
  \langle  g, K_x \rangle_{L^2 (\mu)}=\int_X  g (y) K_y (x) d \mu (y), 
\end{eqnarray*}
where in the last equality we used $K_x (y) =
  \overline{K_y (x)}$. Hence,
\begin{eqnarray*}
  \|  g \|_{L^2 (\mu)}^2 \{ \mu \}_1 & \geq & \langle \tmop{Id}^{\ast}
   g, \tmop{Id}^{\ast}  g \rangle_{H_K}=\langle  g, \tmop{Id}^{\ast}  g \rangle_{L^2 (\mu)}\\
  & = & \int_X \int_X  g (x) \overline{ g (y)} K_y (x) d \mu (x) d
  \mu (y),
\end{eqnarray*}
showing that (i) and (ii) are equivalent.

Testing (ii) on real $ g$'s we have that (iv) holds with $\{ \mu \}_2 =
\{ \mu \}_1$. If viceversa (iv) holds and $ g =  g_R + i  g_I$
is the decomposition of $ g$ in real and imaginary part, then
\begin{equation*}
\begin{split}
  \| I^{\ast}  g \|_{H_K}^2 & \leq (\| I^{\ast}  g_R \|_{H_K}
  + \| I^{\ast}  g_I \|_{H_K})^2 \leq 2 (\| I^{\ast}  g_R
  \|_{H_K}^2 + \| I^{\ast}  g_I \|_{H_K}^2)\\
  & \leq 2 \{ \mu \}_2 \|  g \|_{L^2 (\mu)}^2,
  \end{split}
\end{equation*}
and we obtain the dual form of (i) with $\{ \mu \}_1 \leq 2 \{ \mu \}_2$.
A similar reasoning works if we assume (iv) to hold for positive $ g$'s,
but with a different constant.

It is clear that (vi) and (vii) are equivalent, and (vi) implies (ii) with $\{
\mu \}_1 \leq \{ \mu \}_3$,
\[ \langle \tmop{Id}^{\ast}  g, \tmop{Id}^{\ast}  g \rangle_{H_K} = \langle  g,
   \tmop{Id}^{\ast}  g \rangle_{L^2 (\mu)} \leq \{ \mu \}_3 \|  g
   \|_{L^2 (\mu)}^2 . \]
In the other direction, set $T  g (x) = \int_X K_x (y)  g (y) d \mu
(y) = \tmop{Id}^{\ast}  g (x)$, where $T$ is defined on $L^2 (\mu)$, $T =
T^{\ast}$ and $T$ is positive, $\langle T  g,  g \rangle_{L^2 (\mu)} \geq
0$. Then it has a positive square root $\sqrt{T}$ and by (ii) $\left\|
\sqrt{T}  g \right\|_{L^2 (\mu)}^2 = \langle T  g,  g \rangle_{L^2 (\mu)}
\leq \{ \mu \}_1 \|  g \|_{L^2 (\mu)}^2$, hence (vii) holds with $\{
\mu \}_3 = \{ \mu \}_1$.

Since $\tmop{Re} K$ is positive definite, hence a reproducing kernel, (iii)
and (v) are equivalent for the same reason (i) and (ii) are. Finally, (iii) clearly
implies (iv) and, on the other hand, if (iv) holds, then
\begin{equation*}\begin{split}
  \int_X \int_X  g (x) \overline{ g (y)} \tmop{Re} K_y (x) d \mu (x)
  d \mu (y) & =  \int_X \int_X  g_R (x)  g_R (y) \tmop{Re} K_y (x) d \mu (x) d \mu (y) \\
  & \, \, \, + \int_X \int_X g_I (x)
   g_I (y) \tmop{Re} K_y (x) d \mu (x) d \mu (y) \\ 
  & \leq  \{ \mu \}_2 [\|  g_R \|_{L^2 (\mu)}^2 + \|  g_I
  \|_{L^2 (\mu)}^2] \\
  &= \{ \mu \}_2 \|  g \|_{L^2 (\mu)}^2.
  \end{split}
\end{equation*}
\end{proof}

This lemma allows one to use methods from singular integral theory (where
$\tmop{Re} K$ is the kernel of the ``singular'' integral) on nonhomogeneous
spaces.

This point of view, at least in dyadic theory, started in \cite{ARS2008} in
connection to the problem of Carleson measures for the Drury-Arveson space,
which we will mention below, then taken up by Tchoundia \cite{Tchoundja2008} and Volberg-Wick \cite{VW} in
order to study more general function spaces on the unit ball.

What we aim here, instead, is to give an interpretation of the conformal invariant inequality \eqref{eq:ConformalHardy} in terms of reproducing kernels. Indeed, we will provide such an interpretation for the even more general inequality \eqref{(a)} introduced in Problem \ref{prob:conf inv strong}. The idea is to read \eqref{(a)} as the imbedding inequality of an appropriate RKHS into $L^2( \mu \otimes \mu)$, where $\mu$ is a Borel measure on the rooted tree $T$. Lemma \ref{lem:soft} would then provide various ways to reformulate inequality \eqref{(a)}.

Let $T$ be the rooted tree and $\pi$ an arbitrary positive weight function. We first introduce the \textit{Dirichlet space} $\mathcal{D}
(\pi) := \{ F =\mathcal{I}f : f \in \ell^2 (\pi) \}$, endowed with the inner product $\langle F,G \rangle_{\mathcal{D}
(\pi)}=\langle f,g \rangle_{\ell^2
(\pi)}=\langle \nabla F,\nabla G \rangle_{\ell^2
(\pi)}$, for $F=\mathcal{I}f, G=\mathcal{I}g\in \mathcal{D}
(\pi)$. We claim that this (semi-)Hilbert space has a reproducing kernel,
\begin{equation*}
    K_x (y) = d_{\pi} (x
\wedge y) \assign \sum_{\alpha\in[o^{\ast}, x \wedge y]} \pi (\alpha)^{- 1}.
\end{equation*}
Indeed, for $f\in \ell^2(\pi)$ and $F=\mathcal{I}f$ we have
\begin{equation*}
    F (x)= \sum_{\alpha \in [o^{\ast}, x]} f (\alpha) \pi^{- 1} (\alpha) \pi
  (\alpha)=\sum_{\alpha \in E} f (\alpha) \chi_{[o^{\ast}, x]}(\alpha)\pi^{- 1} (\alpha) \pi
  (\alpha)=\langle F, \chi_{[o^{\ast}, x]}\pi^{- 1} \rangle_{\mathcal{D} (\pi)},
\end{equation*}
from which follows,
\begin{equation*}
    K_x (y) = \mathcal{I}\chi_{[o^{\ast}, x]}\pi^{- 1}(y)=\sum_{\alpha \in [o^{\ast}, y]} \pi^{- 1} (\alpha) \chi
(\alpha \in [o^{\ast}, x]) = d_{\pi} (x \wedge y).
\end{equation*}
It is then imprecise but harmless to say that $\mathcal{D}(\pi)$ is a RKHS. 
The inequality we are re-intepreting as a RKHS imbedding is \eqref{(a)} for $\lambda=\mu \otimes \mu$, i.e.,
\begin{equation}\tag{a'}\label{a'}
   \quad \int_T \int_T | F (x) - F (y) |^2 d \mu(x) d \mu(y) \leq \{
   \mu \} \sum_{\alpha\in E} | \nabla F (\alpha) |^2, 
\end{equation}
We are not quite done yet, since this inequality bounds the $L^2(\mu \otimes \mu)$ norm of the differences of a function with the ($\pi=1$) Dirichlet (semi-)norm of the function itself. However, we argue here that, given a RKHS $H_K$ on a set $X$ such that $1 \in H_K$ (assume $\| 1 \|_{H_K} = 1$), there
is a canonical way to construct the RKHS of its differences, having as
elements the functions $(x, y) \mapsto F (x) - F (y) \backassign  \nabla F (x,
y)$, with $F \in H_K$. We define the kernel $\kappa
: (X \times X) \times (X \times X) \rightarrow \mathbb{C}$ as
\begin{equation}\label{kernel differences}
\begin{split}
  \kappa_{(a, b)} (x, y) & = \kappa ((x, y), (a, b)) \assign \langle K_a - K_b,
  K_x - K_y \rangle\\
  & = K (x, a) - K (y, a) - K (x, b) + K (y, b)\\
  & =  \nabla K_a (x, y) -  \nabla K_b (x, y)=\nabla (K_a-K_b)(x,y).
\end{split}
\end{equation}
We show that $\kappa$ reproduces the space $H_{ \nabla}$ of the functions $ \nabla F$, endowed with
the inner product
\begin{equation}\label{inner product differences}
    \langle  \nabla F,  \nabla G \rangle_{H_{ \nabla}} \assign \langle F - 1 \langle F, 1 \rangle_{H_K}, G - 1 \langle
   G, 1 \rangle_{H_K} \rangle_{H_K},
\end{equation}
which is well defined since $ \nabla F = 0$ if and only if $F$ is constant, i.e., if
and only if $F - 1 \langle F, 1 \rangle_{H_K} = 0$.

\begin{lemma}\label{lem: kernel differences}
  $H_{ \nabla}$ with the inner product \eqref{inner product differences} is a RKHS with kernel $\kappa$ given by \eqref{kernel differences}.
\end{lemma}
\begin{proof}
Since $\langle K_a - K_b, 1 \rangle_{H_K} = 0$, we have
\begin{eqnarray*}
  \langle  \nabla F, \kappa_{(a, b)} \rangle_{H_{ \nabla}} & = & \langle F - 1 \langle F, 1 \rangle_{H_K}, K_a
  - K_b \rangle_{H_K}\\
  & = & F (a) - F (b)\\
  & = &  \nabla F (a, b) .
\end{eqnarray*}
\end{proof}

Lemma \ref{lem: kernel differences} in particular tells us that the space $\mathcal{D}_\nabla$ of differences of functions in the Dirichlet space $\mathcal{D}=\mathcal{D}(1)$ has reproducing kernel
\begin{equation*}
    k_{(a, b)} (x, y) = d (x \wedge a) - d (x \wedge b) - d (y \wedge a) + d (y
   \wedge b) .
\end{equation*}
Inequality \eqref{a'} represents then the boundedness of the imbedding $\mathcal{D}_\nabla\to L^2(\mu \otimes \mu)$, which by means of Lemma \ref{lem:soft} admits various re-writings.

A picture shows that the definition of the kernel of $\mathcal{D}_\nabla$ is independent of the choice of the root,
as we know a priori by conformal invariance.

For  many classical spaces of holomorphic functions, as far as it concerns their imbedding properties, one can substitute the reproducing kernel with its absolute vale causing no losses. It is a natural question if the same applies here.

\begin{problem}
Is \eqref{a'} equivalent to the imbedding in $L^2(\mu \otimes \mu)$ of the space having $|k_{(a, b)} (x, y)|$ as kernel? 
\end{problem}

As a comment to the above problem, we observe that the kernel $k_{(a, b)} (x, y)$ seems to present important cancellations, which might be an indication that tools from
singular integral theory are needed in the characterization of the conformally invariant Hardy's inequality. Indeed, it is a simple exercise to check that for $(a, b), (x, y) \in T \times T$ and $[p, q] := [a, b] \cap [x, y]$, it holds
  \[  k_{(a, b)} (x, y) = \left\{\begin{array}{l}
       + d (p, q) \text{ if $a$ and $x$ (hence, $b$ and $y$) can be joined in
       $T \setminus [p, q]$, }\\
       - d (p, q)  \text{ if $a$ and $x$ (hence, $b$ and $y$) can not be joined
       in $T \setminus [p, q]$.}
     \end{array}\right. \]

\subsection{Quotient structures}\label{sec:quotient}

Dyadic quotient structures appeared for the first time in \cite{ARS2008}, to the best of our knowledge, to deal with
the problem of the Carleson measures for the Drury--Arveson space. Using the
$T^{\ast} T$ argument outlined in Section \ref{sec:rk}, the problem was shown to be equivalent
to the immersion $\tmop{Id} : H_K \rightarrow L^2 (\mu)$ for a tree and a
kernel which we are going to describe in a special case containing all
essential information.

Consider a $4$-adic, rooted tree $T$, whose vertices $x$ at level $d
(o^{\ast}, x) = n$ might be labelled as $4$-adic rationals $x = 0. t_1 \ldots
t_n$, with $t_j \in \mathbb{Z}/ 4\mathbb{Z}$ and an edge joining the parent
$0. t_1 \ldots t_{n - 1}$ with the child $0. t_1 \ldots t_{n - 1} t_n$. Define
similarly the dyadic tree $U$ and consider the surjective map $\Phi : T
\rightarrow U$ induced by the map $[t]_{\tmop{mod} 4} \mapsto [t]_{\tmop{mod}
2}$, sending digits $0, 2$ to binary digit $0$, and digits $1, 3$ to binary
digit $1$.

The map $\Phi$ is a root-preserving {\tmem{tree epimorphism}}: it is
surjective, and $\Phi (x)$ and $\Phi (y)$ are joined by an edge in $U$ if and
only if $x$ and $y$ are joined by an edge in $T$. In other words, we have
defined a quotient structure $U = T / \Phi$ on $T$.

We define a kernel $K_\mathcal{G}$ on $T$ by,

\[ K_\mathcal{G} (x, y) = \frac{2^{- d ([o^{\ast}], [x] \wedge_\mathcal{G} [y])}}{2^{- 2 d (o^{\ast},
   x \wedge_\mathcal{G} y)}}, \]
which can be proved to be positive definite, hence defining a RKHS $H_{K_\mathcal{G}}$. The wedge $\wedge_\mathcal{G}$ is a modified version of the wedge we have used so far. For the exact definition the reader is referred to \cite{ARS2008}. 

The following theorem is proved in  \cite{ARS2008, Tchoundja2008}.

\begin{theorem}\label{DACarleson}
  The following are equivalent for a measure $\mu \geq 0$ on $T$:
  \begin{itemize}
    \item[(i)] The map $\tmop{Id} : H_{K_\mathcal{G}} \rightarrow L^2 (\mu)$ is bounded.
    
    \item[(ii)] We have both the simple condition $\mu (S (\alpha)) \leq C_0
    2^{- d (e (\alpha), o^{\ast})}$ and the inequality
    \[ \int_{S (\alpha)} \left( \int_{S (\alpha)} K_{\mathcal{G}} (x, y) d \mu (y) \right)^s
       d \mu (x) \leq C_p \mu (S (\alpha)), \]
    for one, or equivalently for all, $1 \leq s < \infty$.
  \end{itemize}
\end{theorem}

It is not clear if one needs to introduce the modified wedge $\wedge_\mathcal{G}$ in order for the above theorem to hold. Thus the following problem remains open.

\begin{problem}
  Is it true that Theorem \ref{DACarleson} remains true if we replace the kernel $K_\mathcal{G}$ with the kernel 
  \[ K(x,y):=\frac{2^{- d ([o^{\ast}], [x] \wedge [y])}}{2^{- 2 d (o^{\ast},
   x \wedge y)}}, \]
   where $\wedge$ is the standard tree wedge ?
   
\end{problem}

The real part of
the reproducing kernel of the Drury-Arveson space can be naturally written
down as the quotient of two kernels, which reflect this stratification.
Passing to dyadic decompositions, this leads to the kernels $K_\mathcal{G}$ and $K$ we have just
described, and the Carleson measure problem for the Drury--Arveson space can be
reduced to the theorem stated above.

We have seen that conditions similar to those in the theorem also provide
alternative characterizations of the measures $\mu$ satisfying the Hardy
inequality, at least when $p = 2$. We think that there are here some
interesting questions for further investigation.

\begin{problem}
  Is it possible to have a characterization of the Carleson measures for
  $H_K$ in terms of the potential theory associated with the kernel $K$?
\end{problem}

\subsection{Product structures: poly-trees}\label{sec:product}

The dyadic tree $T$ parametrizes the set of the dyadic subintervals of $[0,
1]$, and the corresponding product structure $T^d$ is defined to parametrize
Cartesian products $R = I_1 \times \ldots \times I_d$ of such intervals:
dyadic rectangles for $d = 2$, etcetera. Following the same lines of Section \ref{subsec:Some potential theory}, a potential theory can now be defined on $T^d$
by taking tensor products of everything on sight, as we will detail below. This leads to a natural extension of the Hardy's inequality to the multi-parameter setting. In this situation, however, characterizing trace measures is a much more complicated problem. We remark that the poly-tree is not a tree, but a graph presenting cycles, and this creates new and major difficulties. So far, solutions to the problem are known for $\sigma\equiv 1$, $p=2$ and for dimension $d =2, 3$ only \cite{AMPS2018}, \cite{AHMV2019}, \cite{MPVZ2020}. It is also known \cite{mozolyako2021multi} that the techniques used in these works are not feasible to be extended to $d=4$ and $p\neq 2$. Let us briefly expand on that.

We identify $T$ with its
vertex set, $T^d \ni x = (x_1, \ldots, x_d)$, and denote by $E^d$ the edge set of $T^d$, $E^d\ni\alpha=(\alpha_1, \ldots, \alpha_d)$. Let $k:\overline{T}\times E\to \mathbb{R}_+$ be the kernel defined in Section \ref{subsec:Some potential theory}.
We define the kernel $\mathbf{k} : \overline{T}^d \times E^d \rightarrow
\{ 0, 1 \}$,
\[ \mathbf{k} (x, \alpha) = \chi_{\lbrace \alpha\supset x\rbrace} (x,\alpha) = \Pi_{j = 1}^d k (x_j, \alpha_j) . \]
Let $\sigma:E^d\to\mathbb{R}_+$ be a positive weight. For a function $\varphi : E^d \rightarrow \mathbb{R}_+$ we set
$\mathbf{I}_\sigma= \mathcal{I}_\sigma \otimes \ldots \otimes \mathcal{I}_\sigma$, i.e.,
\[ \mathbf{I}_\sigma\varphi (x) =  \sum_{\alpha \in E^d} \mathbf{k} (x, \alpha) \varphi (\alpha)\sigma(\alpha)=\sum_{E^d\ni\alpha \supset x} \varphi (\alpha)\sigma(\alpha),
\]
and for $\mu \geqslant 0$ on $\overline{T}^d$,
\[ \mathbf{I}^{\ast} \mu (\alpha) \assign \mu (S(\alpha)), \quad \alpha \in E^d . \]
The $d$-parameter, weighted Hardy inequality for such product structure\footnote{The prototype of the biparameter Hardy inequalities is Sawyer's result
\cite{S1984}, where he considers,
in much more generality, inequalities of the form
\[ \sum_{m, n = 0}^{\infty} \left( \sum_{i = 0}^m \sum_{j = 0}^n f (i, j)
   \right)^2 \mu (m, n) \leqslant [\mu]_{\tmop{Saw}} \sum_{m, n = 0}^{\infty}
   | f (m, n) |^2 . \]
His very clever proof does not extend to the tri-parameter case. Moreover, his inequality is not dyadic, and covers  the facts
here surveyed only in the case of the trivial homogeneous rooted tree $\mathbb{N}$.} is
\[ \int_{\overline{T}^d} \big(\mathbf{I}_\sigma \varphi (x)\big)^p d \mu (x) \leqslant [\mu]_\otimes \sum_{\alpha \in
   E^d} \varphi (\alpha)^p\sigma(\alpha), \quad \varphi \geqslant 0, \]
and the problem is characterizing $\mu$'s for which $[\mu]_\otimes < \infty$, or even
better some geometric, sharp estimate of $[\mu]_\otimes$.

Once we have the kernel $\mathbf{k}$, the general theory \cite{Adams_book} provides us
also with definitions of potentials and energies of measures, and set capacities of
compacta $K \subseteq \overline{T}^d$, as exposed in Section \ref{subsec:Some potential theory}. We can hope at this point that the capacitary estimate does the
job,
\[ \mu \big(\bigcup_{j = 1}^n \alpha^{(j)}\big) \leqslant [\mu]_{\otimes,c} \mathbf{Cap}_{\pi,p} \big(\bigcup_{j = 1}^n \alpha^{(j)}\big), \quad \text{for all } \alpha^{(1)},\ldots,\alpha^{(n)}\in E^d.
\]
Following
Maz'ya's lead, this would follow from a (multi-parameter) {\tmem{Strong Capacitary Inequality}} (see Section \ref{sec:SCI}),
\begin{equation}\tag{SCI}
    \sum_{k = - \infty}^{+ \infty} 2^{2 k} \mathbf{Cap}_{\pi,p} (x: \mathbf{I}f (x) > 2^k) \leqslant A \| \varphi
   \|_{\ell^p_+ (E^d,\pi)}^p .
\end{equation}
Here a major difficulty appears: the standard proofs of (SCI) depend, more or
less explicitly, on the {\tmem{boundedness principle}} for potentials of
measures,
\[ \sup \{ \mathbf{V}_p^{\mu,\sigma} (x) : x \in \overline{T}^d \} \leqslant B \cdot
   \max \{ \mathbf{V}_p^{\mu,\sigma} (x) : x \in \tmop{supp} (\mu) \}, \]
but in the multi-parameter situation such principle miserably fails.

\begin{proposition}[\cite{AMPS2018}]
  For $d \geqslant 2$ there exist measures $\mu^K$ which are equilibrium for a
  compact $K \subset (\partial T)^d$, hence automatically satisfy
  $\mathbf{V}^{\mu^K}:=\mathbf{V}_2^{\mu^K,1} \leqslant 1$ on $\tmop{supp} (\mu^K)$, such that
  $\mathbf{V}^{\mu^K} (x) = + \infty$ at some point $x \in (\partial T)^d$. 
\end{proposition}

The idea is to have a set $K$ which is rarefied, but ``curved'' is such a way
many ``not too thin'' rectangles join it to the point $x$, like rays
focusing on it.

The way out of this difficulty, implemented in \cite{AMPS2018} for $d =
2$, $p=2$ and $\sigma\equiv 1$, is proving a {\tmem{distributional boundedness principle}}. 

\begin{theorem}
  There is $C > 0$ such that for $\lambda > 1$ and for an equilibrium measure
  $\mu$,
  \[ \mathbf{Cap} (\{ x : \mathbf{V}^{\mu} (x) > \lambda \}) \leqslant C
     \frac{\| \mathbf{I}^{\ast} \mu \|_{\ell^2}^2}{\lambda^{2 + 1}} . \]
\end{theorem}

The inequality would follow by rescaling if one has $2$ instead of $2 + 1$.
This weaker form of the boundedness principle suffices to produce a variation
of a classical proof of the Strong Capacitary Inequality,
\[ [\mu]_\otimes \approx [\mu]_{\otimes,c}, \quad \text{for } d = 2. \]
This result was then extended to $d = 3$, but no higher, in \cite{MPVZ2020}. With some major difficulty, the capacitary characterization of the measures
for which the multiparameter, dyadic Hardy inequality holds is true at least
for $d = 2, 3$. What about the other characterizations and proofs?

It is proved in \cite{AHMV2019} that a mass--energy condition holds as well in $d =
2$, and in \cite{MPVZ2020} this was extended to $d = 3$, always for $p=2$. More precisely, for $d =
2, 3$ we have that $[\mu]_\otimes \approx [[\mu]]_\otimes$, where $[[\mu]]_{\otimes}$ is the
best constant in
\begin{equation}\tag{ME$\otimes$}\label{eq:ME multiparameter}
    \sum_{E^d\ni\beta \subseteq \alpha} \mu (S(\beta))^2 \leq [[\mu]]_\otimes \  \mu (S(\alpha))
   < \infty, \quad \alpha\in E^d.
\end{equation}

This fact might surprise practitioners of the Hardy space on the bidisc. It
was proved in \cite{C1974} that Carleson measures for the Hardy space on the
bidisc are not characterized by a ``single-box condition'' such as \eqref{eq:ME multiparameter}, and
A. Chang proved in \cite{Chang1979} that the characterization holds if one allows
multiple boxes. One might expect that a multiple box condition like
\[ \sum_{\beta \subseteq \cup_{j = 1}^n \alpha^{(j)}} \mu (S(\beta))^2 \leq [\mu]_{\tmop{mult}}
   \mu \big(\bigcup_{j = 1}^n S(\alpha^{(j)})\big) < \infty, \quad \text{for all } \alpha^{(1)},\ldots,\alpha^{(n)}\in E^d, \]
might not be weakened, but if fact this is not the case.

The proofs we surveyed for the one-parameter Hardy operator seem not to work
in the multi-parameter case. The simple maximal proof, for instance, does not
work because, contrary to the usual dyadic, weighted maximal function, its
several parameter versions,
\[ \mathcal{M}_{\mu} f (x) = \sup_{E^d\ni\alpha \supset x} \frac{1}{\mu (S(\alpha))}
   \int_{S(\alpha)} f d \mu, \]
are not necessarily weakly bounded on $L^1$, neither they are bounded on
$L^2$. In the unweighted case, the $L^2$ boundedness of the multiparameter
maximal function was proved in \cite{JMZ1935}, and a nice account of multiparameter
theory with applications to martingales and the Hardy space is in \cite{G1983}.

\begin{problem}
It would be interesting to know whether, like in the one parameter case, for $1 \leqslant s < \infty$,
\[ \sup_{\alpha\in E^d}\frac{1}{\mu (S(\alpha))} \int_{S(\alpha)} \left( \int_{S(\alpha)} \delta_d (x \wedge y) d \mu (y)
   \right)^s d \mu (x) \approx C_s[[\mu]]^s , \]
where, $\delta_d (x \wedge y) \assign \prod_{j=1}^dd (x_j \wedge y_j)$ .
\end{problem}

\section{Appendix} 

\subsection{Dyadic decomposition of Ahlfors regular metric spaces}\label{sec:dyadic}


Let $\mathcal{D}$ be the sets of the dyadic intervals $I_{n, j} = [(j - 1) /
2^n, j / 2^n)$ in $[0, 1)$, $1 \leq j \leq 2^n$, $n \geq 0$.  Let also $\mu, \nu : \mathcal{D} \to \mathbb{R}_+$ be two weight functions. Then we say that the \textit{two-weight dyadic Hardy's inequality} holds, if there exists a positive constant $[\mu]$, possibly depending on $\pi$ and $p$, such that for all $f : \mathcal{D} \to \mathbb{R}_+$,  

\begin{equation} \label{HardyDyadicIntervals} \sum_{I \in \mathcal{D}} \left( \sum_{J \supseteq I} f
  (J) \right)^p \mu (I) \leq [\mu] \sum_{I \in \mathcal{D}} f
  (I)^p \nu (I).  \end{equation}

It is clear that the above is, in fact, a Hardy's inequality on the homogeneous dyadic tree: interpret $\mathcal{D}$ as the vertex set of a tree, where two vertices are connected by an edge if and only if one of the corresponding intervals $I_{n,j}, I_{m,k}$ contains the other and $|n-m|=1$, and set $o$ be the vertex corresponding to $[0,1)$. 
Observe that here, as compared to the general formulation \eqref{HardyInequality}, we are in the simplest case when $\mu$ is supported on the tree $\mathcal{D}$ rather than on $\overline{\mathcal{D}}$. In the paper we chose to work in higher generality and allow trees to be  not necessarily homogeneous and the measure $\mu$ to give mass also to $\partial T$, the natural boundary of the tree. In this way it becomes clear that the symmetric-space structure of the group of automorphisms of the homogeneous tree plays no role in (most of) this theory. Moreover, to extend the support of the measure up to boundary is justified by the problem of exceptional sets at the boundary, i.e.,  of trace measures.

We remark that $\mathcal{D}$ is just the prototype of decomposition of a metric space. In a more general context, if $ (X,\mu, \rho) $ is a homogeneous metric measure space, then its Christ's decomposition \cite[Theorem 11]{Christ1990} provides a family of generalized ``cubes'' $\{Q_k^\alpha \}$ which can be readily checked to form a tree.

Let us show that \eqref{HardyDyadicIntervals} is a genuine generalization of the classical Hardy inequality \eqref{classicHardy}.
Suppose for example that $\{I_n\}$ is an infinite  branch of dyadic intervals, i.e., $I_n \in \mathcal{D}, I_0=[0,1), I_{n+1}\subseteq I_n $ and $ 2|I_{n+1}|=|I_n|.$ Set also $\pi(I)=\mu(I) = 0 $ if $I\in \mathcal{D}$ is not one of the $I_n$  and $\mu(I_n)=:U_n, \pi(I_n)=:V_n$. Write $\varphi_n$ for $\varphi (I_n)$. Then the dyadic Hardy's inequality takes the form 
\[ \sum_{n=0}^\infty \Big( \sum_{m=0}^n \varphi_m \Big)^p U_n \leq [\mu] \sum_{n=0}^\infty \varphi_n^p V_n, \]
which is of course the discrete analogue of Muckenhoupt's two-weight Hardy's inequality in $\mathbb{R}_+.$ In particular, by choosing $U_n = n^{-p} $ and $V_n=1$ one gets back to \eqref{classicHardy}.

An interesting example of an unexpected application of \eqref{HardyDyadicIntervals}, coming  from complex analysis, is the problem of characterizing Carleson measures for Besov spaces $B^p_a$. A Carleson measure for $B^p_a$ is a positive Borel measure
$\tilde{\mu}$ on the unit disc $\mathbb{D}$ of the complex plane for which there exists
a constant $C(\tilde{\mu}) < \infty$ such that, for all $f$ holomorphic on $\mathbb{D}$,
\[ \tag{Besov}\label{Dirichlet:Carleson} \int_{\mathbb{D}} | f |^p d \tilde{\mu} \leq C(\tilde{\mu}) \left[ | f (0) |^2 + \int_{\mathbb{D}} |  f' (z)
  |^p (1-|z|^2)^{p+a-2} d x d y \right], \]
with $0 \leq ap < 1$. Such problems appear in connection to the
characterization of multipliers and of exceptional sets at the boundary for
spaces of holomorphic functions, sequences of interpolation, and more. The first result is by Stegenga \cite{Stegenga1980} who manages to characterize such measures for $p=2$  in terms of a condition involving  Riesz capacities  of compact subsets of the unit disc. The root of Stegenga's work can be traced back to earlier work of Maz'ya \cite{Mazya1973} and Adams \cite{Adams1976}. For the case $1<p<\infty$, a similar characterization of Carleson measures in terms of non linear Riesz capacities was later obtained by Verbitskii \cite{verb} and rediscovered by Wu \cite{Wu1999}.

More recently it  was proved in \cite{Kerman1988, Arcozzi2002} that the Carleson inequality for Besov spaces is equivalent to a dyadic Hardy's inequality. More precisely, for a given dyadic interval $I = I_{n, j} \in \mathcal{D}$, let $Q
 (I) $ be the set of points $z=re^{i\theta}$ with $\theta/2\pi \in I $ and $ 1-2^{-n}\leq r \leq 1-2^{-n-1} $. Set also $\mu (I) := \tilde{\mu} (Q
 (I))$ and $\pi (I) = 2^{- a n}$. Then $ \tilde{\mu} $ is a Carleson measure for $B^p_a, 0\leq ap<1$, if and only if the triple $p, \mu, \pi$ satisfies the dyadic Hardy's inequality \eqref{HardyDyadicIntervals}. Motivated by this application we shall call Carleson measures (or trace measures) all measures $\mu : \mathcal{D} \to \mathbb{R}_+$ satisfing the Hardy's inequality on trees \eqref{HardyInequality}.




The dyadic setting is much ductile. The same inequality \eqref{HardyDyadicIntervals}, with
different choices of the weight $\pi$, can
be used to characterize Carleson measures for holomorphic spaces in several
dimensions or for spaces of harmonic functions, trace inequalities for
potential spaces, and more. Many such problems, in fact, can be proven to be
equivalent to their dyadic counterparts, and often \eqref{HardyDyadicIntervals} is the form they
assume.

\subsection{Bessel potentials on the boundary of the dyadic tree}\label{subsec:Bessel}

 The content of this section is specific to the homegenous tree. Out of simplicity, we consider only the dyadic case, but everything we say applies, mutatis mutandis, to homogeneous trees of any degree.

Our objective in this section is to introduce, using as usual the axiomatic theory of Adams and Hedberg, a seemingly different potential theory on the boundary of the dyadic tree which depends on two parameters $p$ and $s$. Subsequently we will use the inequality of Muckenhoupt and Wheeden in order to prove that it is ``equivalent'' to the potential theory introduced in Section \ref{subsec:Some potential theory} for the same parameter $p$ and a particular choice of the weight $\pi$. Equivalent means that for compact sets that lie on the boundary of the dyadic tree, the capacities of the sets measured by means of the two different potential theories are comparable.

Let $T$ be the dyadic tree and consider the compact Hausdorff space  $\overline{T}$ and the measure space $(\partial T,dx)$. Fix some parameter $0 < s \leq 1/p^*$. We then define the $s-$Bessel kernel as  \[G_s(x,y):=|\alpha|^{-s}, \quad x,y\in \partial T\]
where $\alpha$ is the (unique) edge such that $e(\alpha)=x\wedge y.$ 

Following, as always, \cite{Adams_book} we define the Bessel potential  of a function $\varphi : \partial T \to \mathbb{R}_+$ as 
\[ G_s\varphi (y) : = \int_{\partial T} G_s(x,y) \varphi(x) dx, \quad y\in \overline{T}.  \]
Also the Bessel co--potential of a measure $\mu $ is defined as 
\[ G^*_s\mu (x) : = \int_{\overline{T}}G_s(x,y)d\mu(y), \quad x\in\partial T. \]

Notice that, if temporarily we use the notation $\hat{\alpha}$ to denote  the only child of $\alpha$ lying in $[o^*,x]$, for some fixed vertex $x$, we can estimate the Bessel co--potential as follows 

\begin{align*}
    G_s^*\mu(x)& =\int_{\overline{T}}G_s(x,y)d\mu(y)=\sum_{\alpha\supset x}\int_{ S(\alpha)\setminus  S(\hat{\alpha})}|\alpha|^{-s}d\mu(y) \\ 
    &=\sum_{\alpha\supset x}(\mu( S(\alpha))- \mu( S(\hat{\alpha})))|\alpha|^{-s}.
\end{align*}
Hence, we trivially have $ G^*_s\mu(x)\leq \sum_{\alpha\supset x}\mu( S(\alpha))|\alpha|^{-s}$. On the other hand,
\begin{align*}
    G_s^* \mu(x) & =  \sum_{\alpha\supset x}\frac{\mu( S(\alpha))}{|\alpha|^s}-\sum_{\alpha\supset x}\frac{\mu( S(\hat{\alpha}))}{2^s|\hat{\alpha}|^s}=\mu(\partial T)+(1-2^{-s})\sum_{\alpha\supset x, \alpha\neq \omega}\frac{\mu( S(\alpha))}{|\alpha|^s} \\ 
   &  \geq (1-2^{-s}) \sum_{a\supset x} \frac{\mu(S(\alpha))}{|\alpha|^s}.
\end{align*}
In other words,  \[ G_s^*\mu(x)\approx\sum_{\alpha\supset x}\mu( S(\alpha))|\alpha|^{-s}. \]
The associated Bessel energy is given by
\begin{align*}
    E_p^s(\mu) & :=\int_{\overline{T}}\Big(\int_{\partial T}G_s(x,y)G^*_s\mu(x)^{p^*-1}dx\Big)d\mu(y) \\ 
    & =
    \int_{\partial T} G_s^*\mu(x)^{p^*}dx  \approx \int_{\partial T}\Big(\sum_{\alpha\supset x}\frac{\mu( S(\alpha))}{|\alpha|^s}\Big)^{p^*}dx  \\
    & \approx \sum_\alpha \frac{\mu( S(\alpha))^{p^*}}{|\alpha|^{sp^*-1}} =  \mathcal{E}_{p,\pi}(\mu),
\end{align*}
where $\pi(\alpha) = | \alpha |^{\frac{1-sp^*}{1-p^*}}$. Notice that we have used the Muckenhoupt--Wheeden inequality \eqref{eq:MW} in the last step.

Therefore, since the energies associated to a positive Borel measure via the two different potential theories are comparable, we can conclude that the $\tmop{Cap}_{p,\pi}$ capacity of a compact subset of the boundary of the dyadic tree and its $s-$Bessel capacity are comparable. In particular compact sets of zero $\tmop{Cap}_{p,\pi}$ capacity coincide with those of zero $s-$Bessel capacity. 




\bibliographystyle{plain}
\bibliography{biblio}

\end{document}